\documentclass[12pt]{article}

\usepackage[utf8]{inputenc}
\usepackage{latexsym,amsfonts,amssymb,amsthm,amsmath}

\usepackage{graphicx}
\usepackage{xcolor}
\usepackage{algorithm}
\usepackage{algpseudocode}
\usepackage{authblk}
\usepackage{hyperref}

\title{Pessimistic bilevel optimization approach for decision--focused learning}
\author[1,2]{Diego Jiménez}
\author[1]{Bernardo K. Pagnoncelli}
\author[2]{Hande Yaman}
\affil[1]{SKEMA Business School}
\affil[2]{Faculty of Economics and Business, KU Leuven}

\begin{document}

\maketitle

\begin{abstract}
    The recent interest in contextual optimization problems, where randomness is associated with side information, has led to two primary strategies for formulation and solution. The first, \textit{estimate-then-optimize}, separates the estimation of the problem’s parameters from the optimization process. The second, \textit{decision-focused optimization}, integrates the optimization problem's structure directly into the prediction procedure. In this work, we propose a pessimistic bilevel approach for solving general decision-focused formulations of combinatorial optimization problems. Our method solves an $\varepsilon$-approximation of the pessimistic bilevel problem using a specialized cut generation algorithm. We benchmark its performance on the 0-1 knapsack problem against estimate-then-optimize and decision-focused methods, including the popular SPO+ approach. Computational experiments highlight the proposed method’s advantages, particularly in reducing out-of-sample regret.
\end{abstract}

\vspace{0.5in}

\section{Introduction}
\label{intro}
In many real-world applications, decisions must be made in the presence of uncertainty, where some key parameters remain unknown at the time of decision making. Typically, measurements of contextual variables, known as covariates or features, are available. While not directly part of the decision model, these covariates are often correlated with the unknown parameters and can be leveraged to produce reliable estimates.

For instance, in vehicle routing problems, determining the least-cost routes for a fleet of vehicles requires anticipating future traffic conditions, which are inherently uncertain. These conditions can be estimated using contextual information such as the current weather, time of day, and the state of the road network \cite{LIU2023, islam2023}. Similarly, in portfolio investment problems, decisions on allocating a limited budget across assets with uncertain returns can benefit from macroeconomic indicators, news sentiment, or even social media trends as sources of contextual information \cite{xu-cohen-2018, pagnoncelli2023synthetic, wang2023, pena2024modified}.

In this setting, contextual optimization (CO) provides a formal mathematical framework for addressing such problems, as outlined in \cite{sadana2023survey}. This approach assumes the existence of a joint probability distribution linking the unknown parameters and the observed features. Similarly to stochastic optimization models, CO aims to minimize the expected cost of the final decision with respect to this distribution, conditioned on the observed feature values:
\begin{equation}
    \mathbf{z}^{CO}(\mathbf{x}) \in \underset{\mathbf{z} \in \mathcal{Z}}{\text{arg min }} \mathbb{E}_{\mathbb{P}(\mathbf{y}|\mathbf{x})}  [ c(\mathbf{y}, \mathbf{z}) ]. \label{eq:main_prob}
\end{equation}
Here, $\mathbf{y} \in \mathcal{Y} \subset \mathbb{R}^{n_y}$ denotes the vector of random parameters, which are correlated with features $\mathbf{x} \in \mathcal{X} \subset \mathbb{R}^{n_x}$. This correlation is described by the conditional probability distribution $\mathbb{P}(\mathbf{y}|\mathbf{x})$. In this formulation, the downstream optimization problem (DOP), for a fixed value of $\mathbf{y}$, corresponds to $\underset{\mathbf{z} \in \mathcal{Z}}{\text{min }} c(\mathbf{y}, \mathbf{z})$, in which $\mathbf{z}$ corresponds to the vector of decision variables and $ \mathcal{Z}$ is the associated feasible set. Note that in formulation \eqref{eq:main_prob}, random parameters appear only in the objective function, which is the case addressed in this study.

Problem \eqref{eq:main_prob} indicates that the optimal CO decision, $\mathbf{z}^{CO}(\mathbf{x})$, referred to as a \textit{predictive prescription} in \cite{Bertsimas2020}, depends on the observed context, represented by the feature vector $\mathbf{x}$, and the conditional probability distribution $\mathbb{P}(\mathbf{y}|\mathbf{x})$. However, in real-world applications, these distributions are typically unavailable. Instead, decision makers rely on historical datasets of features and unknown parameters, denoted as $\mathcal{D} = \lbrace (\mathbf{x}_i, \mathbf{y}_i) \rbrace_{i = 1}^n$. Consequently, DOP problems require approximations, many documented in a recent survey \cite{sadana2023survey}. The most common approach involves creating \textit{estimation models} to predict the unknown parameters, which are then used to solve the CO problem.

When function $c(\mathbf{y}, \mathbf{z})$ is linear in $\mathbf{y}$, the expectation in \eqref{eq:main_prob} simplifies as  $\mathbb{E}_{\mathbb{P}(\mathbf{y}|\mathbf{x})}  [ c(\mathbf{y}, \mathbf{z}) ] = c(\mathbb{E}_{\mathbb{P}(\mathbf{y}|\mathbf{x})}  [\mathbf{y}], \mathbf{z})$.
To approximate the conditional expectation $\mathbb{E}_{\mathbb{P}(\mathbf{y}|\mathbf{x})}[\mathbf{y}]$, an \textit{estimator} $f:\mathcal{X} \rightarrow \mathbb{R}^{n_y}$ is computed,
 selected from a family of functions $f \in \mathcal{F}$. Using this estimator, the approximated CO solution can be expressed as:
\begin{equation}
    \mathbf{z}^\star (f(\mathbf{x})) \in \underset{\mathbf{z} \in \mathcal{Z}}{\text{arg min }}  c(f(\mathbf{x}), \mathbf{z}). \label{eq:main_aprox}
\end{equation}
In the \textit{machine learning} (ML) context, these estimators are calibrated by minimizing a non-negative loss function $\ell: \mathbb{R}^{n_y} \times \mathbb{R}^{n_y} \rightarrow \mathbb{R}+$, which measures the distance between true and predicted values based on the available dataset. The empirical risk minimization (ERM) problem is formulated as:
\begin{alignat}{2} f^\star \in \underset{f \in \mathcal{F}}{\text{arg min }} & \dfrac{1}{n} \sum_{i=1}^n \ell \big ( f(\mathbf{x}_i), \mathbf{y}_i \big ) + \Omega(f), \label{eq:erm}
\end{alignat}
where $\Omega: \mathcal{F} \rightarrow \mathbb{R}+$ is a regularization term designed to prevent overfitting, such as the lasso penalty. A common choice for $\ell$ in regression problems is the mean squared error (MSE), defined as $\ell_{MSE} (f(\mathbf{x}_i), \mathbf{y}_i) = || f(\mathbf{x}_i) - \mathbf{y}_i ||^2$.

A wide range of ML predictors benefits from well-established and computationally efficient algorithms to solve \eqref{eq:erm} for $\ell_{MSE}$. However, this approximation to CO problems, referred to as \textit{estimate-then-optimize}, overlooks  the existence of the subsequent decision model problem
on the estimator's computation.
As a result, the optimal estimator \( f^\star \) is decoupled from \( c(\mathbf{y}, \mathbf{z}) \) and \( \mathcal{Z} \), even though those factors can heavily influence the choice of  \( \mathbf{z}^\star \big( f^\star (\mathbf{x}) \big) \) and, by extension, the final decision cost \( c \big( \mathbf{y}, \mathbf{z}^\star \big( f^\star (\mathbf{x}) \big) \big) \).

In this context, a growing body of research has focused on calibrating estimator models to minimize final decision errors. This methodology is referred to in the literature as \textit{integrated estimation-optimization} (IEO) by \cite{elmachtoub2023estimate}, \textit{decision-focused learning} by  \cite{mandi2024}, and \textit{problem-aware} loss functions by \cite{Malaga2022}. By embedding the optimization structure of the decision problem into the loss function $ \ell(\cdot)$ definition, this approach accounts for decision errors and produces \emph{decision-focused} estimators.

The most intuitive way to embed the structure of the DOP into the loss function is to define it directly in terms of the regret relative to the optimal solution with full information \(\mathbf{y}\):
\begin{equation}
        \Tilde{\ell}_{IEO} ( f(\mathbf{x}), \mathbf{y}) := c \big (\mathbf{y}, \mathbf{z}^\star ( f(\mathbf{x}) ) \big) - \underset{\mathbf{z} \in \mathcal{Z}}{\text{min }} c \big (\mathbf{y}, \mathbf{z} \big ),  \label{eq:ieo_loss_amb}
\end{equation}
where \(\mathbf{z}^\star \big ( f(\mathbf{x}) \big )\) is the optimal solution of the DOP using the estimation \(f(\mathbf{x})\) as the cost vector. In \eqref{eq:ieo_loss_amb}, estimation errors are evaluated through their impact on the final decision cost, effectively integrating both predictive and prescriptive components into the loss. However, solving the empirical risk minimization problem \eqref{eq:erm} with the loss function \eqref{eq:ieo_loss_amb} often results in highly challenging optimization problems, even for relatively simple DOPs.

In \cite{ban2019big} the authors propose  ERM and kernel-based methods to address a contextual version of the newsvendor problem. They demonstrate that failure to incorporate features into the problem results in suboptimal decisions that are not asymptotically optimal as the number of data points approaches infinity.
Similarly,  \cite{Elmachtoub2022} highlights that estimation-focused and decision-focused estimators for the same datasets can differ significantly. The authors argue that minimizing decision errors directly requires integrating the estimation and optimization models into a unified framework. However, this integration transforms \eqref{eq:erm} into a more complex mathematical programming problem, which can hinder the practical application of the technique. Specifically, when no analytical expression of $\mathbf{z}^\star \big ( f(\mathbf{x}) \big )$ exists for a given DOP, solving the ERM with the loss function $\Tilde{\ell}_{IEO}(\cdot)$ leads to a bilevel problem, as \eqref{eq:ieo_loss_amb} inherently embeds the DOP minimization.

Note that definition \eqref{eq:ieo_loss_amb} can be ambiguous if the minimizer $\mathbf{z}^{\star} (f(\mathbf{x}))$ is not unique. Within the bilevel optimization framework, selecting a minimizer $\mathbf{z}^{\star} (f(\mathbf{x}))$ that also minimizes the leader's objective function, in this case, the ERM cost function, is referred to as the \textit{optimistic} case. Conversely, in the \textit{pessimistic} case, the worst-case minimizer $\mathbf{z}^{\star} (f(\mathbf{x}))$ is chosen, as determined by the highest leader's objective cost.

The case of continuous optimization is examined in \cite{Malaga2022}, where the lower-level problem is reformulated using Karush-Kuhn-Tucker (KKT) conditions, resulting in a single-level formulation for the optimistic case of the loss function \eqref{eq:ieo_loss_amb}. However, as the optimistic case can produce degenerate solutions, a quadratic regularization term is incorporated into the objective function to address this issue.

For the linear programming (LP) case, \cite{Elmachtoub2022} and \cite{Bucarey2024} propose approximations to develop tractable estimators for the problem. Both studies address the pessimistic version of \eqref{eq:ieo_loss_amb}, which results in an even more challenging problem, necessitating the use of approximations. In \cite{Elmachtoub2022}, a surrogate function of the pessimistic IEO loss, known as the SPO+ (smart predict-then-optimize) loss, is employed. This approach reduces the original problem to a linear single-level formulation by dualizing the internal DOP. Conversely, \cite{Bucarey2024} does not rely on a surrogate loss but instead represents the lower-level problem using primal and dual feasibility conditions alongside strong duality. However, the resulting formulation is a non-convex quadratic problem, for which approximations and heuristics are applied.

Asymptotic guarantees are examined in \cite{Fatma2022}, which investigates the conditions under which prediction methods ensure strong out-of-sample performance. Specifically, the study considers both MSE and SPO+ loss functions, representing optimization-agnostic and IEO methods, respectively. The results highlight that the absence of Fisher consistency, which depends both on the distribution $\mathbb{P}(\mathbf{y}|\mathbf{x})$ and the loss function $\ell(\cdot)$, can lead to poor out-of-sample performance. Interestingly, while Fisher consistency is achieved by MSE under certain conditions, it is not guaranteed for SPO+. However,  SPO+ still outperforms MSE in specific scenarios, such as when the DOP corresponds to the  continuous knapsack problem.
In \cite{Estes2023}, the authors focus on two-stage continuous linear CO problems, where uncertain parameters appear exclusively on the right-hand side of the second-stage constraints. They demonstrate that the SPO+ loss is not applicable in this setting and introduce a novel convex IEO loss function, for which sub-gradients can be computed.

The application of IEO to combinatorial optimization problems remains underexplored in the literature. One notable exception is \cite{mandi2019smart}, where the SPO+ loss is applied to the linear relaxations of combinatorial problems.
Our review highlights two significant gaps in the field. First, no exact method has been proposed for solving the ERM of the pessimistic case of $\eqref{eq:ieo_loss_amb}$; current approaches rely solely on approximations. Second, combinatorial optimization problems involving integer variables have received limited attention, with a notable lack of exact methods for addressing this case.

To address these gaps, we identify the following key contributions of this work:\begin{enumerate}
    \item We propose a method for producing estimators by minimizing the true pessimistic IEO loss without relying on approximations or surrogate functions. This approach also allows for minimizing the IEO loss in purely combinatorial optimization problems without requiring a convex hull description, unlike SPO+.
    \item To solve the resulting model, we adapt a column-and-constraint generation approach that iteratively solves an integer master problem and a series of integer subproblems. To enhance computational performance, we introduce a novel branch-and-cut algorithm that solves a single master problem. Computational comparisons reveal that the branch-and-cut approach significantly outperforms the column-and-constraint generation method.
    \item Numerical experiments are conducted to evaluate the out-of-sample performance of the proposed method, benchmarking it against linear regression (LR) and SPO+. The results demonstrate that for problems of moderate dimensions, our approach outperforms SPO+ unless the convex hull is explicitly computed. For larger problems, our method consistently outperforms both LR and SPO+.
\end{enumerate}

\section{Mathematical framework and formulations} \label{sec:math_form}
\subsection{Preliminaries and notation} \label{ssec:prelim}
The following notation is used. Vectors are represented by bold lowercase characters, matrices by uppercase characters, and scalars by lowercase characters. Specifically, $\mathbf{x}$, $\mathbf{y}$, and $\mathbf{z}$  represent the vectors of features, random parameters, and decision variables, respectively. Estimators, denoted as $f_\mathbf{w}$, are  uniquely defined by a parameter vector $\mathbf{w} \in \mathbb{R}^{n_w}$. Thus, optimizing over $f_\mathbf{w}$ within the family of functions $\mathcal{F}$ is equivalent to optimizing over  $\mathbf{w}$ in $\mathbb{R}^{n_w}$, so we focus on optimizing over $\mathbf{w}$ instead of $f_\mathbf{w}$. Uppercase calligraphic characters denote sets, and for any $n \in \mathbb{N}$, $[n] = \lbrace 1, 2, ..., n \rbrace$. Additionally, $\mathcal{Z}$  represents the feasible set of the DOP, and $\mathcal{Z}^\star (\mathbf{y})$ the set of minimizers of the problem $\underset{\mathbf{z} \in \mathcal{Z}}{\text{min }} c(\mathbf{y}, \mathbf{z})$.

\subsection{Bilevel integrated estimation-optimization} \label{ssec:bilevel_ieo}
Here, we analyze the ERM minimization of the IEO loss \eqref{eq:ieo_loss_amb}. As pointed out in the previous section, the definition of $\Tilde{\ell}_{IEO}(\cdot)$ is ambiguous when $|\mathcal{Z}^\star (f_{\mathbf{w}} (\mathbf{x}))|
\allowbreak > 1$ and it requires a selection rule. In the context of bilevel optimization, optimistic and pessimistic selection rules represent the two extreme and the most studied cases.

When considering the optimistic case, the ERM \eqref{eq:erm} can lead to degenerate solutions. As mentioned in \cite{Malaga2022} and  \cite{Elmachtoub2022}, optimistic IEO solutions may produce policies mapping different features values to the same decisions. In particular, if $\Omega(f_\mathbf{w}) = 0$ in \eqref{eq:erm}, i.e., no regularization is considered, and if there exists some $\Tilde{\mathbf{w}}$ for which $f_{\Tilde{\mathbf{w}}} (\mathbf{x}) = 0$, i.e., the null vector for all $\mathbf{x} \in \mathcal{X}$, then we get that $\mathcal{Z}^\star (f_{\Tilde{\mathbf{w}}} (\mathbf{x})) = \mathcal{Z}$, resulting in an IEO loss equal to $0$. However, such a predictor would most likely lead to very low-quality out-of-sample results.

In \cite{Malaga2022}, a regularization term $\Omega(f_\mathbf{w})$ is added to avoid this case, increasing the complexity of the original problem. In \cite{Elmachtoub2022, Bucarey2024}, as in our work, the pessimistic approach is taken. With this approach, the associated loss function reads:
\begin{equation}
    \ell_{IEO} \big ( f_{\mathbf{w}}(\mathbf{x}), \mathbf{y} \big ) = \underset{\mathbf{z} \in \mathcal{Z}^\star (f_{\mathbf{w}}(\mathbf{x}))}{\text{max }} c \big (\mathbf{y}, \mathbf{z} \big )  -  \underset{\mathbf{z} \in \mathcal{Z}}{\text{min }} c \big (\mathbf{y}, \mathbf{z} \big ). \label{eq:ieo_loss}
\end{equation}
This choice of the pessimistic approach serves a dual purpose. First, it avoids the pathological case discussed earlier without introducing additional terms to the loss function. Second, it can be interpreted as a worst-case optimization approach. However, pessimistic bilevel problems are notoriously challenging to solve and, in certain instances, may even be unsolvable, as highlighted in \cite{Wiesemann2013}.

\subsection{Pessimistic bilevel IEO formulation} \label{ssec:pess_b_ieo}
In our setting, as DOP, we consider a pure linear combinatorial optimization problem (COP) of the form:
\begin{subequations}\label{eq:bMILP}
\begin{alignat}{3}
    \underset{\mathbf{z}}{\text{min }} &
    \mathbf{y}^\intercal \mathbf{z} &\label{eq:bMILP1} \\
    \text{s.t. } & \mathbf{A}\mathbf{z} \geq \mathbf{b}, \label{eq:bMILP2} \\
    & \mathbf{z} \in \lbrace 0, 1 \rbrace^{n_y}.  \label{eq:bMILP3}
\end{alignat}
\end{subequations}
As in \cite{Elmachtoub2022}, we consider the case where the cost vector $\mathbf{y}$ is the only random parameter of the optimization model. The information available about it is contained in a dataset $\mathcal{D} := \lbrace (\mathbf{x}_i , \mathbf{y}_i) \rbrace_{i=1}^n$, in which $\mathbf{x}_i$ and $\mathbf{y}_i$ represent the $i$-\textit{th} historical observations of features and cost vectors samples, respectively.
The main goal of the present proposal is to fit a predictor $f_{\mathbf{w}}(\mathbf{x})$, by minimizing the average normalized pessimistic IEO loss \eqref{eq:ieo_loss} over $\mathcal{D}$. In that regard, our main problem to solve is the following:
\begin{alignat}{3}
    \mathbf{w}^\star \in \underset{\mathbf{w}}{\text{arg min }} &  \Bigg \lbrace \dfrac{1}{n} \sum_{i = 1}^n  \Bigg [\underset{\mathbf{z}_i \in \mathcal{Z}^\star (f_{\mathbf{w}}(\mathbf{x}_i))}{\text{max }} \lbrace \bar{\mathbf{y}}^{\intercal}_{i} \mathbf{z}_i \rbrace - 1 \Bigg ] \Bigg \rbrace.  \label{eq:spo_min_full}
\end{alignat}
where $\bar{\mathbf{y}}_i = \frac{\mathbf{y}_i}{\underset{\mathbf{z} \in \mathcal{Z}}{\text{min }} \mathbf{y}_{i}^{\intercal} \mathbf{z}}$ is the normalized cost vector.
To simplify the presentation, the constant term 1 will be omitted from \eqref{eq:spo_min_full} in subsequent formulations, without affecting the set of optimal solutions.

Problem \eqref{eq:spo_min_full} expresses that $\mathbf{w}^{\star}$ is the predictor parameter vector that minimizes the average worst-case decision regret that results from using $f_{\mathbf{w}} (\mathbf{x})$ as an estimation of the real cost vector $\mathbf{y}$ in solving the DOP over the dataset $\mathcal{D}$. Note that this interpretation is in line with Theorem 3 in \cite{Goerigk2024}, which states that, under certain conditions, an instance of a worst-case regret optimization problem can be solved by addressing an instance of a pessimistic bilevel optimization.

Unlike other IEO approaches such as \cite{Malaga2022}, \cite{Elmachtoub2022}, and \cite{Bucarey2024}, we employ a normalized IEO loss for computing the ERM. This choice is motivated by the fact that the normalized loss is commonly used to evaluate out-of-sample performance, as demonstrated in \cite{Elmachtoub2022}. Thus, using the normalized IEO loss for ERM helps to provide a clearer understanding of both in-sample and out-of-sample behavior.

Formulation \eqref{eq:spo_min_full} presents three main challenges. First, pessimistic bilevel problems are, in general, hard to solve \cite{Wiesemann2013}, and in our case we have a three-level nested structure, which becomes more visible when \eqref{eq:spo_min_full} is rewritten as
\begin{subequations}\label{eq:SPO_pess}
\begin{alignat}{2}
\underset{\mathbf{w}}{\text{min }} & \underset{\lbrace \mathbf{z}_i \rbrace_{i=1}^n}{\text{max }}
\dfrac{1}{n} \sum_{i = 1}^n \bar{\mathbf{y}}_{i}^{\intercal} \mathbf{z}_i   \label{eq:SPO_pess1}\\
\text{s.t. } &  \mathbf{z}_i \in \mathcal{Z}_{i}^{\star} (f_{\mathbf{w}} (\mathbf{x}_i)) = \underset{\mathbf{z}}{\text{arg min }} \Big \lbrace f_{\mathbf{w}}(\mathbf{x}_i)^\intercal \mathbf{z} \, : \, \mathbf{z} \in \mathcal{Z} \Big \rbrace, \quad \forall i \in [n], \label{eq:SPO_pess2} \\
& \mathbf{w} \in \mathbb{R}^{n_w}, \label{eq:SPO_pess3}
\end{alignat}
\end{subequations}
where $\mathcal{Z} = \big \lbrace \mathbf{z} \in \lbrace 0,1 \rbrace^{n_y} , \; \mathbf{A} \mathbf{z} \geq \mathbf{b} \big \rbrace$. Here, $\mathbf{w}$ can be interpreted as the vector of the leader's decision variables, and $\mathbf{z}_i$, for $i \in [n]$, as follower's variables. Still, it is worth mentioning that the bilevel structure of \eqref{eq:SPO_pess} is not inherited by its context or application, as in interdiction problems or Stackelberg games, in which leaders and followers represent real agents' interactions. Rather, the bilevel structure of \eqref{eq:SPO_pess} is a consequence of the use of the loss \eqref{eq:ieo_loss}.

Second, follower's variables, i.e., $\lbrace \mathbf{z}_i \rbrace_{i=1}^n$, are binary, which prohibit the use of any method relying on continuous functions \cite{Dempe2018}. Finally, the follower's objective function $f_{\mathbf{w}} (\mathbf{x}_i)^\intercal \mathbf{z}$ is non-linear (bilinear in the best case, if $f_{\mathbf{w}} (\mathbf{x}_i)$ is linear in $\mathbf{w}$).

Upon reviewing the available methods in the literature for solving bilevel mixed integer linear programming (MILP) models, it is evident that, as noted in \cite{KLEINERT2021}, most exact approaches are based on the so-called \textit{high-point relaxation} (HPR). This technique optimizes over the set $\mathcal{H}$, comprising variables that satisfy both leader and follower constraints \cite{Koppe2010, Fischetti2017, Tahernejad2020}.
In the optimistic IEO case, this set corresponds to:
\begin{equation} \label{eq:hpr_set}
    \mathcal{H} = \Big \lbrace (\mathbf{w}, \lbrace \mathbf{z}_i \rbrace_{i=1}^n) \in \mathbb{R}^{n_w} \times   \lbrace 0,1 \rbrace^{n_y}  \times ...  \times  \lbrace 0,1 \rbrace^{n_y} : \,  \mathbf{A} \mathbf{z}_i \geq \mathbf{b}, \forall i \in [n] \Big \rbrace.
\end{equation}
Note that if we optimize $\frac{1}{n} \sum_{i=1}^n \bar{\mathbf{y}}_{i}^{\intercal} \mathbf{z}_i$ over $\mathcal{H}$, the predictor parameters $\mathbf{w}$ could take any  value. In fact, the HPR destroys the structure of \eqref{eq:SPO_pess} since $\mathbf{w}$ is only present on the follower's objective function.

Alternatively, an iterative solution scheme is presented in \cite{Wiesemann2013} for solving convergent $\varepsilon$-approximations of general pessimistic bilevel problems. Such a scheme could be applied in principle to solve \eqref{eq:SPO_pess}, without the shortcomings of HPR. However, it would require solving highly non-linear optimization problems consecutively, rendering its application unpractical, even for small values of $n$.
More recently, \cite{BoZeng2020} presented valuable results for general pessimistic bilevel problems and proposed a column-and-constraint generation algorithm for the mixed integer linear case. However, due to the non-linearity of the follower's cost function, $f_{\mathbf{w}} (\mathbf{x})^\intercal \mathbf{z}$, their methodology and some of their results cannot be directly applied to our case. Consequently, we build upon their findings for the general case to develop a solution method for \eqref{eq:SPO_pess}.

\section{Solution methodology}
\label{sec:solmeth}
\subsection{Valid relaxations for pessimistic IEO} \label{ssec:valid_relax}
In \cite{BoZeng2020}, the author proposes the following relaxation:
\begin{subequations}\label{eq:bzV2}
\begin{alignat}{3}
\underset{}{\text{min }} & \dfrac{1}{n} \sum_{i=1}^n \theta_i & \label{eq:bzV2_1} \\
\text{s.t. }  & \mathbf{w} \in \mathbb{R}^{n_w}, \label{eq:bzV2_2} \\
& \bar{\mathbf{z}}_i \in \mathcal{Z} \quad & \forall i \in [n],  \label{eq:bzV2_3} \\
& \bar{\mathcal{Z}}_i(\mathbf{w}, \bar{\mathbf{z}}_i) =  \Big \lbrace  \mathbf{z} \in \mathcal{Z} \, : \, f_{\mathbf{w}} (\mathbf{x}_i)^\intercal (\mathbf{z} - \bar{\mathbf{z}}_i) \leq 0  \Big \rbrace  \quad & \forall i \in [n], \label{eq:bzV2_5} \\
& \theta_i \geq \text{max } \Big \lbrace  \bar{\mathbf{y}}_{i}^\intercal \mathbf{z}_i \, : \, \mathbf{z}_i \in \bar{\mathcal{Z}}_i(\mathbf{w}, \bar{\mathbf{z}}_i) \Big \rbrace \quad & \forall i \in [n],  \label{eq:bzV2_4}
\end{alignat}
\end{subequations}
whose main advantage is to replace \eqref{eq:SPO_pess2}  by a more tractable constraint $\mathbf{z}_i \in \bar{\mathcal{Z}}_i(\mathbf{w}, \bar{\mathbf{z}}_i)$ and proves it to be a tight relaxation  of \eqref{eq:SPO_pess}.
To give an idea of the argument in \cite{BoZeng2020} behind the relaxation, observe that for any feasible point $(\mathbf{w}^\star, \lbrace \mathbf{z}_{i}^\star \rbrace_{i=1}^n)$ of \eqref{eq:SPO_pess}, if we select $\bar{\mathbf{z}}_{i}^\star = \mathbf{z}_{i}^\star$ for all $i \in [n]$, we get that $\bar{\mathcal{Z}}_i (\mathbf{w}^\star, \bar{\mathbf{z}}_{i}^{\star}) = \mathcal{Z}_{i}^{\star} (f_{\mathbf{w}^\star} (\mathbf{x}_i)) $. Hence, $\mathbf{z}_{i}^{\star} \in \underset{\mathbf{z}}{\text{arg max }}  \lbrace \bar{\mathbf{y}}_{i}^{\intercal} \mathbf{z} \, : \, \mathbf{z} \in \bar{\mathcal{Z}}_i (\mathbf{w}^{\star}, \bar{\mathbf{z}}_{i}^{\star}) \rbrace $, proving that any feasible solution of \eqref{eq:SPO_pess} can be extended to one of $\eqref{eq:bzV2}$. In addition, for any other $\bar{\mathbf{z}}_{i}\in \mathcal{Z}$, we have $ \mathcal{Z}_{i}^{\star} (f_{\mathbf{w}^\star} (\mathbf{x}_i)) \subseteq \bar{\mathcal{Z}}_i (\mathbf{w}^{\star}, \bar{\mathbf{z}}_i)$. The relaxation is tight since for any $\mathbf{w}^\star$, there exists a $\bar{\mathbf{z}}^{\star}_{i}$, optimal for $\mathbf{w}^\star$ in the sense of \eqref{eq:bzV2}, such that $\text{max}\lbrace \bar{\mathbf{y}}_{i}^{\intercal} \mathbf{z} \, : \, \mathbf{z} \in \bar{\mathcal{Z}}_i (\mathbf{w}^\star, \bar{\mathbf{z}}_{i}^{\star}) \rbrace = \text{max}\lbrace \bar{\mathbf{y}}_{i}^{\intercal} \mathbf{z} \, : \, \mathbf{z} \in \mathcal{Z}_{i}^{\star} (f_{\mathbf{w}^\star} (\mathbf{x}_i)) \rbrace$ for all $i \in [n]$.
\cite{BoZeng2020} provided solution algorithms for mixed integer linear pessimistic bilevel models.  However, no solution method was provided for the case in which the follower's objective function depends on the leader's variables. In fact, in the current configuration, problem \eqref{eq:bzV2} involves a non-linear cost in the follower's objective function, represented by the product $f_{\mathbf{w}} (\mathbf{x}_i)^\intercal \mathbf{z}$.
For this reason, an extension should be made to cover this case, allowing the resolution of the ERM of the pessimistic IEO loss.


Note that if $\mathbf{z}_i$ can only take integer values and $\mathcal{Z}$ is bounded, then $\mathcal{Z}$ has a finite number of points. Let $\mathcal{F}(\mathcal{Z}) = \lbrace \mathbf{z}^e \rbrace_{e=1}^{n_p}$ be the full enumeration set of the feasible points of $\mathcal{Z}$. Consequently, inequalities \eqref{eq:bzV2_5}-\eqref{eq:bzV2_4} can be rewritten as:
\begin{equation} \label{eq:ep_rep}
    \theta_i \geq \text{max } \Big \lbrace  \bar{\mathbf{y}}_{i}^\intercal \mathbf{z}^e \, : \, \mathbf{z}^e \in \mathcal{F}(\mathcal{Z}), f_{\mathbf{w}} (\mathbf{x}_i)^\intercal (\mathbf{z}^e - \bar{\mathbf{z}}_i) \leq 0  \Big \rbrace, \quad \forall i \in [n].
\end{equation}

Our goal is to exploit this discrete representation's capabilities to design an iterative method. To this end, consider $\mathcal{Z}_{i}^{(k)} \subseteq \mathcal{F}(\mathcal{Z})$ for $i \in [n]$, and also the following sets:
\begin{align*}
    \Xi_{i}^{(k)}(\mathbf{w}, \bar{\mathbf{z}}_i) =& \Big \lbrace \mathbf{z}^{e} \in \mathcal{Z}_{i}^{(k)}, \, f_{\mathbf{w}} (\mathbf{x}_i)^\intercal (\mathbf{z}^e - \bar{\mathbf{z}}_i) \leq 0 \Big  \rbrace, \quad \forall i \in [n] \\
    \Gamma^{(k)} =& \Big \lbrace (\mathbf{w}, \lbrace \theta_i \rbrace_{i=1}^n, \lbrace \bar{\mathbf{z}}_i \rbrace_{i=1}^n) \, : \, \eqref{eq:bzV2_2}-\eqref{eq:bzV2_3}, \\
    & \quad\quad
    \theta_i \geq \text{max } \big \lbrace  \bar{\mathbf{y}}_{i}^\intercal \mathbf{z}^e \, : \, \mathbf{z}^e \in \Xi_{i}^{(k)}(\mathbf{w}, \bar{\mathbf{z}}_i) \big \rbrace, \; \forall i \in [n]  \Big \rbrace.
\end{align*}
With these definitions in hand, we can state the following result:

\subsubsection*{Lemma 1: Valid relaxation of Pessimistic IEO} \label{lemma1}
Consider the problem $P^{(k)}: \text{min} \lbrace \frac{1}{n} \sum_{i=1}^n \theta_i \, ; \, (\mathbf{w}, \lbrace \theta_i \rbrace_{i=1}^n, \lbrace \bar{\mathbf{z}}_i \rbrace_{i=1}^n) \in \Gamma^{(k)} \rbrace$. If $\mathcal{Z}_{i}^{(k)} \subseteq \mathcal{Z}_{i}^{(k+1)}$ for all $i \in [n]$, then $P^{(k)}$ is a relaxation of $P^{(k+1)}$.\\

\begin{proof}
As both $P^{(k)}$ and $P^{(k+1)}$ have the same objective function, it suffices to check whether $\Gamma^{(k+1)} \subseteq \Gamma^{(k)}$. Moreover, as equations \eqref{eq:bzV2_2}-\eqref{eq:bzV2_3} are the same for both sets, we only need to verify that if $(\mathbf{w}, \lbrace \theta_i \rbrace_{i=1}^n, \lbrace \bar{\mathbf{z}}_i \rbrace_{i=1}^n) \in \Gamma^{(k+1)}$, then the inequalities $\theta_i \geq \text{max} \big \lbrace  \bar{\mathbf{y}}_{i}^\intercal \mathbf{z}^e \, : \, \mathbf{z}^e \in \Xi_{i}^{(k)} (\mathbf{w}, \bar{\mathbf{z}}_i)  \big \rbrace$ hold for all $i \in [n]$. Since  $\Xi_{i}^{(k)}(\mathbf{w}, \bar{\mathbf{z}}_i)\subseteq \Xi_{i}^{(k+1)}(\mathbf{w}, \bar{\mathbf{z}}_i)$, we have $ \text{max } \big \lbrace  \bar{\mathbf{y}}_{i}^\intercal \mathbf{z}^e \, : \, \mathbf{z}^e \in \Xi_{i}^{(k+1)}(\mathbf{w}, \bar{\mathbf{z}}_i) \big \rbrace \geq  \text{max } \big \lbrace  \bar{\mathbf{y}}_{i}^\intercal \mathbf{z}^e \, : \, \mathbf{z}^e \in \Xi_{i}^{(k)}(\mathbf{w}, \bar{\mathbf{z}}_i) \big \rbrace$. As $\theta_i\geq  \text{max } \big \lbrace  \bar{\mathbf{y}}_{i}^\intercal \mathbf{z}^e \, : \, \mathbf{z}^e \in \Xi_{i}^{(k+1)}(\mathbf{w}, \bar{\mathbf{z}}_i) \big \rbrace$, we have that $\theta_i\geq  \text{max } \big \lbrace  \bar{\mathbf{y}}_{i}^\intercal \mathbf{z}^e \, : \, \mathbf{z}^e \in \Xi_{i}^{(k)}(\mathbf{w}, \bar{\mathbf{z}}_i) \big \rbrace$ by $\Xi_{i}^{(k)}(\mathbf{w}, \bar{\mathbf{z}}_i)\subseteq \Xi_{i}^{(k+1)}(\mathbf{w}, \bar{\mathbf{z}}_i)$, proving the result.
\end{proof}

The result obtained in Lemma 1 allows us to construct a framework in which sequentially tighter relaxations can be obtained by the successive addition of feasible points
of $\mathcal{Z}$ to $\mathcal{Z}_{i}^{(k)}$ for each $i \in [n]$, resulting in a sequence of non-decreasing lower bounds. Additionally, upper bounds are easy to compute, since, given any $\mathbf{w} \in \mathbb{R}^{n_w}$, a feasible pessimistic bilevel solution of \eqref{eq:bzV2} can be calculated by solving the following problem for each $i \in [n]$:
\begin{subequations}\label{eq:spo_sp}
\begin{alignat}{3}
SP_i (\mathbf{w}):
\, \mathbf{z}_{i}^{SP} (\mathbf{w}) \in\underset{\mathbf{z}_{i}}{\text{arg max }} &
\bar{\mathbf{y}}_{i}^{\intercal} \mathbf{z}_i  &  \label{eq:spo_sp1}\\
\text{s.t. } & f_{\mathbf{w}}(\mathbf{x}_i)^\intercal \mathbf{z}_i \leq \text{opt} \big ( f_{\mathbf{w}}(\mathbf{x}_i) \big ), \label{eq:spo_sp2} \\
&  \mathbf{z}_i \in \mathcal{Z},  \label{eq:spo_sp3}
\end{alignat}
\end{subequations}
where $\text{opt} \big ( f_{\mathbf{w}}(\mathbf{x}_i) \big ) = \text{min} \lbrace f_{\mathbf{w}}(\mathbf{x}_i)^\intercal \mathbf{z}_i  : \, \mathbf{z}_i \in \mathcal{Z} \rbrace$. Then, an upper bound can be computed as $\frac{1}{n} \sum_{i=1}^n \bar{\mathbf{y}}_{i}^{\intercal} \mathbf{z}_{i}^{SP}(\mathbf{w})$.

In the next section, we describe in detail the basis of our framework that allows to solve \eqref{eq:SPO_pess} to optimality.

\subsection{Column-and-constraint generation algorithm for pessimistic IEO} \label{ssec:ccg}
The next step is to explicitly write the previously proposed relaxation $P_k : \text{min} \lbrace \frac{1}{n} \sum_{i=1}^n \theta_i  ; \allowbreak\, (\mathbf{w}, \lbrace \theta_i \rbrace_{i=1}^n, \lbrace \bar{\mathbf{z}}_i \rbrace_{i=1}^n) \in \Gamma^{(k)} \rbrace$. First, notice that, for each $i \in [n]$, constraint $\theta_i \geq \text{max} \big \lbrace  \bar{\mathbf{y}}_{i}^\intercal \mathbf{z}^e \, : \, \mathbf{z}^e \in \Xi_{i}^{(k)} (\mathbf{w}, \bar{\mathbf{z}}_i)  \big \rbrace$ could be rewritten as
\begin{equation}\label{eq:disj}
    \theta_i \geq \bar{\mathbf{y}}_{i}^{\intercal} \mathbf{z}^{e}  - M \mathbb{I} \big [ f_{\mathbf{w}}(\mathbf{x}_i)^\intercal (\mathbf{z}^{e} - \bar{\mathbf{z}}_i)  > 0 \big ] \quad \forall \mathbf{z}^{e} \in \mathcal{Z}_{i}^{(k)},
\end{equation}
with $M$ as a big-$M$ value, and $\mathbb{I}[\cdot]$ an indicator function. To be able to implement that disjunction by making use of off-the-shelf solvers, we propose the following model:
\begin{subequations}\label{eq:bzV3}
\begin{alignat}{3}
P^{(k)}_{\varepsilon}: \underset{}{\text{min }} &
\dfrac{1}{n} \sum_{i=1}^n \theta_i & \label{eq:bzV3_1} \\
\text{s.t. }
& \theta_i \geq \bar{\mathbf{y}}_{i}^{\intercal} \mathbf{z}^{e} \lambda_{ie} &  \forall e \in \mathcal{Z}_{i}^{(k)}, \, \forall i \in [n], \label{eq:bzV3_2} \\
& f_{\mathbf{w}}(\mathbf{x}_i)^\intercal (\mathbf{z}^{e} - \bar{\mathbf{z}}_i)  \geq  \varepsilon - M\lambda_{ie}  & \quad \forall e \in \mathcal{Z}_{i}^{(k)}, \, \forall i \in [n], \label{eq:bzV3_4} \\
& \lambda_{ie} \in \lbrace 0,1 \rbrace & \forall e \in \mathcal{Z}_{i}^{(k)}, \; \forall i \in [n],  \label{eq:bzV3_5} \\
& \bar{\mathbf{z}}_i \in \mathcal{Z} & \forall i \in [n],
\label{eq:bzV3_6} \\
& \theta_i \geq \theta_{i}^{min}  & \forall i \in [n], \label{eq:bzV3_7} \\
& \mathbf{w} \in \mathbb{R}^{n_w}. \label{eq:bzV3_8}
\end{alignat}
\end{subequations}

The value of the binary variable $\lambda_{ie}$ indicates whether $\mathbf{z}^e \in \Xi_{i}^{(k)} (\mathbf{w}, \bar{\mathbf{z}}_i)$, or not, by setting its value to 1 or 0, respectively. This condition is modeled by constraints \eqref{eq:bzV3_4}, which are checking whether the cost of $\mathbf{z}^{e}$ is not larger than the cost of $\bar{\mathbf{z}}_i$, with respect to the predicted cost vector $f_{\mathbf{w}}(\mathbf{x}_i)$. If $\lambda_{ie}=0$, then constraints \eqref{eq:bzV3_4} become $f_{\mathbf{w}}(\mathbf{x}_i)^\intercal (\mathbf{z}^{e} - \bar{\mathbf{z}}_i)  \geq  \varepsilon$, imposing the cost of $\mathbf{z}^{e} $ to be strictly larger than the one of $\bar{\mathbf{z}}_i$ (the value of $\varepsilon$ should be chosen accordingly). If $\lambda_{ie}=1$, then this constraint becomes redundant ($M$ can be set equal to $\varepsilon$ or larger)  and constraint \eqref{eq:bzV3_2} is activated. 
In summary, constraints \eqref{eq:bzV3_2} and \eqref{eq:bzV3_4} ensure that  $\theta_i \geq \text{max} \big \lbrace  \bar{\mathbf{y}}_{i}^\intercal \mathbf{z}^e \, : \, \mathbf{z}^e \in \Xi_{i}^{(k)} (\mathbf{w}, \bar{\mathbf{z}}_i)  \big \rbrace$ holds for all $i \in [n]$.
Additionally, the lower bound values $\theta_{i}^{min}$ can be obtained by computing $\theta_{i}^{min} = \underset{\mathbf{z} \in \mathcal{Z}}{\text{min }} \bar{\mathbf{y}}_{i}^{\intercal} \mathbf{z}$.
Note that we can add the valid inequalities $\theta_i \geq \bar{\mathbf{y}}_{i}^{\intercal} \bar{\mathbf{z}}_i$  for all $i \in [n]$ to the formulation because of constraints
\eqref{eq:bzV2_4}.


Thereafter, as formulation \eqref{eq:bzV3} is an $\varepsilon$--approximation of a  valid relaxation of \eqref{eq:bzV2}, which is, at the same time, a valid relaxation of the original problem \eqref{eq:SPO_pess}, we can construct the following general column-and-constraint generation (CCG) algorithm for the pessimistic IEO problem, similar to the one proposed in \cite{BoZeng2020}:

\begin{algorithm}[h!]\label{alg:genCCG0}
\caption{A general CCG algorithm for Pessimistic IEO}
\vskip6pt
\begin{algorithmic}[1]\label{alg:genCCG}
\State \textbf{Input}: Dataset $ \mathcal{D} = \lbrace(\mathbf{x}_i , \mathbf{y}_i) \rbrace_{i=1}^n$.  Initial $\mathcal{Z}_{i}^{(0)}, \, \forall i \in [n]$. Algorithm tolerance $\Delta > 0$. Set $k = 0$.
    \While{$ UB - LB \geq \Delta$}
        \State Solve $P^{(k)}_{\varepsilon}$ to get $\theta^{\star} $ and $\mathbf{w}^\star$
        \State Update $LB \leftarrow \frac{1}{n} \sum_{i=1}^n \theta_{i}^{\star}$
        \State Solve $SP_{i}(\mathbf{w}^\star)$ and get $\mathbf{z}_{i}^{SP}(\mathbf{w}^\star)$ for each $i \in [n]$
        \State Update $\mathcal{Z}_{i}^{(k+1)}\leftarrow \mathcal{Z}_{i}^{(k)} \cup \lbrace \mathbf{z}_{i}^{SP} (\mathbf{w}^{\star}) \rbrace$ for each $i \in [n]$
        \State Update $UB \leftarrow \frac{1}{n}\sum_{i=1}^n \bar{\mathbf{y}}_{i}^{\intercal} \mathbf{z}_{i}^{SP}(\mathbf{w}^\star) $
        \State Set $k \leftarrow k+1$

    \EndWhile
    \State \textbf{Output}: Predictor parameters vector $\mathbf{w}^\star$.
\Statex
\end{algorithmic}
\end{algorithm}

What is left is to prove that the Algorithm 1 converges after a finite number of iterations, which is done by the following result.

\subsubsection*{Lemma 2: Convergence of the CCG algorithm}
Algorithm 1 converges, after a finite number of iterations, to a pessimistic IEO solution $\mathbf{w}^\star$ whose gap from the optimal is at most $\Delta$.\\

\begin{proof}
     The only case in which Algorithm 1 could iterate indefinitely is if  the subproblems $SP_{i} (\mathbf{w}^\star)$ compute solutions $\mathbf{z}_{i}^{SP} (\mathbf{w}^\star)$ that are already contained in $\mathcal{Z}_{i}^{(k)}$ for all $i \in [n]$, and yet the gap $UB - LB$ is still greater than or equal to the predefined tolerance $\Delta$.
    Assume that $\mathbf{z}_{i}^{SP} (\mathbf{w}^\star)\in \mathcal{Z}_{i}^{(k)}$ for all $i \in [n]$. Let $i\in [n]$. Since  $\bar{\mathbf{y}}_{i}^{\intercal} \mathbf{z}_{i}^{SP} (\mathbf{w}^\star) = \text{max} \lbrace  \bar{\mathbf{y}}_{i}^\intercal \mathbf{z}^e \, : \, \, \mathbf{z}^e \in \mathcal{Z}_{i}^{\star}(\mathbf{w}^\star) \rbrace$, we have that $f_{\mathbf{w}} (\mathbf{x}_i)^\intercal (\mathbf{z}_{i}^{SP} (\mathbf{w}^\star)- \bar{\mathbf{z}}_{i}^{\star}) \leq 0$. This, combined  with     $\mathbf{z}_{i}^{SP} (\mathbf{w}^\star) \in \mathcal{Z}_{i}^{(k)}$, implies that
    $\mathbf{z}_{i}^{SP} (\mathbf{w}^\star) \in \Xi_{i}^{(k)} (\mathbf{w}^\star, \bar{\mathbf{z}}_{i}^{\star})$. As $\theta_{i}^\star \geq   \text{max}  \lbrace  \bar{\mathbf{y}}_{i}^\intercal \mathbf{z}^e \, : \, \mathbf{z}^e \in \Xi_{i}^{(k)} (\mathbf{w}^\star, \bar{\mathbf{z}}_{i}^{\star})  \rbrace $, we have $  \theta_{i}^\star \geq \bar{\mathbf{y}}_{i}^{\intercal} \mathbf{z}_{i}^{SP} (\mathbf{w}^\star)$.
    Summing up these inequalities over all $i \in [n]$ and dividing by $n$, we obtain $\frac{1}{n}\sum_{i=1}^n  \theta_{i}^\star \geq \frac{1}{n}\sum_{i=1}^n \bar{\mathbf{y}}_{i}^{\intercal} \mathbf{z}_{i}^{SP}(\mathbf{w}^\star)$. The left-hand side represents LB and the right-hand side represents UB. This demonstrates that if   $\mathbf{z}_{i}^{SP} (\mathbf{w}^\star)\in \mathcal{Z}_{i}^{(k)}$ for all $i \in [n]$, then $LB\geq UB$, satisfying the stopping condition.
\end{proof}

\subsection{Branch--and--cut algorithm for pessimistic IEO} \label{ssec:bac_gen}

One potential drawback of using Algorithm 1 is that the master problem \eqref{eq:bzV3},  a mixed integer program, is solved at every iteration. This can result in significant computational effort. At each iteration, the CCG algorithm introduces new variables $\lambda_{ie}$, and constraints \eqref{eq:bzV3_2}--\eqref{eq:bzV3_4} to the master problem, thereby strengthening the relaxation \eqref{eq:bzV3}. However, these constraints can be incorporated in a more efficient manner.

Motivated by this observation, we propose a branch--and--cut scheme to solve \eqref{eq:bzV3}, with $\mathcal{Z}_{i}^{(k)} = \mathcal{F}(\mathcal{Z})$ for all $ i \in [n]$, corresponding to the pessimistic IEO problem. In this approach, constraints \eqref{eq:bzV3_2}--\eqref{eq:bzV3_4} are added lazily. Unlike Algorithm 1, which adds both variables and constraints iteratively, the branch--and--cut scheme defines all $\lambda_{ie}$ variables during the initialization phase, without associating any $\mathbf{z}^{e}$ to them. As the algorithm progresses, points $\mathbf{z}^{e}$ are computed dynamically, and the corresponding constraints \eqref{eq:bzV3_2}--\eqref{eq:bzV3_4} are appended accordingly. If the constraints associated with a variable $\lambda_{ie}$ have not yet been added  for specific values of $i$ and $e$, that variable does not influence the model solution. Thus, whenever an integer solution is encountered in the branch--and--bound tree, the constraints \eqref{eq:bzV3_2}--\eqref{eq:bzV3_4}, which are valid global cuts, can be added. Moreover, even when fractional solutions are found, global valid cuts can be added by solving \eqref{eq:spo_sp} for the current value of $\mathbf{w}$.

In the present implementation, we use the branch--and--cut framework provided by an off-the-shelf solver. 
Therefore, branching and node selection stages are managed by the solver, while the contribution of the present proposal is on the cut generation procedure, which we describe below.

\begin{enumerate}
    \item \textbf{Initialization}: Define the set of initial points $\mathcal{Z}_{i}^{(0)} \subset \mathcal{F}(\mathcal{Z})$ for all $ i \in [n]$.  Set $e(i) = |\mathcal{Z}_{i}^{(0)}|$  for all $ i \in [n]$, the number of initial added points. Then, define the initial problem:
    \begin{alignat*}{3}
    \underset{}{\text{min }} &
    \dfrac{1}{n} \sum_{i=1}^n \theta_i &  \\
    \text{s.t. }
    & \theta_i \geq \bar{\mathbf{y}}_{i}^{\intercal} \mathbf{z}^{e} \lambda_{ie} &  \forall e \in \mathcal{Z}_{i}^{(0)}, \, \forall i \in [n],  \\
    & f_{\mathbf{w}}(\mathbf{x}_i)^\intercal (\mathbf{z}^{e} - \bar{\mathbf{z}}_i)  \geq  \varepsilon - M\lambda_{ie}  & \quad \forall e \in \mathcal{Z}_{i}^{(0)}, \, \forall i \in [n], \\
    & \lambda_{ie} \in \lbrace 0,1 \rbrace & \forall e \in \mathcal{F}(\mathcal{Z}), \; \forall i \in [n],  \\
    & \eqref{eq:bzV3_6}, \eqref{eq:bzV3_7}, \eqref{eq:bzV3_8}.
    \end{alignat*}
    \item \textbf{Cut generation}: Whenever a solution is obtained  by solving the linear relaxation at a node (either fractional or integer) with a lower objective function value than the current global upper bound, obtain the current solution $(\mathbf{w}, \lbrace \theta_{i} \rbrace_{i \in [n]}, \lbrace \bar{\mathbf{z}}_{i} \rbrace_{i \in [n]}, \lbrace \lambda_{ie} \rbrace_{i \in [n], e \in \mathcal{F}(\mathcal{Z})})$. Compute $\mathbf{z}_{i}^{SP} (\mathbf{w})$ from \eqref{eq:spo_sp}, for all $i \in [n]$. If $\mathbf{z}_{i}^{SP} (\mathbf{w}) \in \mathcal{Z}^{(k)}_i$ for some $i \in [n]$, then set $\mathcal{Z}_{i}^{(k+1)} \leftarrow \mathcal{Z}_{i}^{(k)}$, and continue. Else, update $ \mathcal{Z}_{i}^{(k+1)} \leftarrow  \mathcal{Z}_{i}^{(k)} \cup \lbrace \mathbf{z}_{i}^{SP} (\mathbf{w}) \rbrace$, $e(i) \leftarrow e(i) + 1$, and add the following constraints:
        \begin{subequations} \label{eq:bac_cuts}
        \begin{align}
        & \theta_{i} \geq \bar{\mathbf{y}}_{i}^{\intercal} \mathbf{z}_{i}^{SP} (\mathbf{w}) \lambda_{i \; e(i)}, \label{eq:bac_cuts_1}  \\
        & f_{\mathbf{w}}(\mathbf{x}_i)^\intercal (\mathbf{z}_{i}^{SP} (\mathbf{w}) - \bar{\mathbf{z}}_{i})  \geq  \varepsilon - M\lambda_{i \; e(i)}.  \label{eq:bac_cuts_2}
        \end{align}
        \end{subequations}
    Then, update the current global upper bound if $\frac{1}{n} \sum_{i=1}^n \bar{\mathbf{y}}_{i}^{\intercal} \mathbf{z}_{i}^{SP} (\mathbf{w})$ is lower than the current one. Otherwise, continue.
\end{enumerate}

Notice that cuts \eqref{eq:bac_cuts_1}--\eqref{eq:bac_cuts_2} are exactly the same as the constraints \eqref{eq:bzV3_2}--\eqref{eq:bzV3_4}, respectively. 
Note that, even when getting fractional solutions,  points of $\mathcal{Z}$, $\mathbf{z}_{i}^{SP} (\mathbf{w}) \in \mathcal{F}(\mathcal{Z})$, can be generated by solving \eqref{eq:spo_sp}, for any value of $\mathbf{w}$, and, therefore, pessimistic IEO feasible solutions, with $\theta_i (\mathbf{w}) = \mathbf{y}_{i}^{\intercal} \mathbf{z}_{i}^{SP} (\mathbf{w})$, and $\lambda_{ie} = 1$ for $\mathbf{z}_{ie} = \mathbf{z}_{i}^{SP}(\mathbf{w})$, and $\lambda_{ie} = 0$ otherwise, for all $i \in [n]$.  

\subsection{Additional settings} \label{ssec:settings}

In this section, some practicalities when solving the pessimistic IEO problem are discussed.

\subsubsection{Linearization of $f_{\mathbf{w}}(\mathbf{x}_i)^{\intercal} \bar{\mathbf{z}}_i$} First of all, it is important to remark that no assumptions have been specified on the predictor function $f_{\mathbf{w}}(\mathbf{x})$, which is present in constraints \eqref{eq:bzV3_4}. As $\bar{\mathbf{z}}_i$ is binary, it is possible to write down an equivalent reformulation of \eqref{eq:bzV3}, linear in $f_{\mathbf{w}}(\mathbf{x})$ and $\bar{\mathbf{z}}_i$, if bounds on $f_{\mathbf{w}}(\mathbf{x})$ are imposed. Particularly, if $f_{min} \leq f_{\mathbf{w}}(\mathbf{x})_j \leq f_{max}$, we define a new continuous variable $\mathbf{u}_i \in \mathbb{R}^{n_y}$ such that
\begin{subequations}\label{eq:mc_relax}
\begin{align}
& f_{min} \bar{z}_{ij} \leq u_{ij} \leq f_{max} \bar{z}_{ij} & \forall i \in [n], j \in [n_y], \\
& f_{min} (1-\bar{z}_{ij}) \leq f_{\mathbf{w}}(\mathbf{x})_j - u_{ij} \leq  f_{max} (1-\bar{z}_{ij}) & \forall i \in [n], j \in [n_y].
\end{align}
\end{subequations}
In this case, it can be ensured that $u_{ij} = f_{\mathbf{w}} (\mathbf{x})_j \bar{z}_{ij}$ when getting an integer solution. Finally, if $f_{\mathbf{w}}(\mathbf{x})$ is linear in $\mathbf{w}$, we can solve the pessimistic IEO by solving a MILP formulation.

\subsubsection{Warmstart initialization} Points of $\mathcal{Z}$ can be obtained beforehand, and included at the initialization stage of the algorithm to accelerate its resolution. In particular, by repeatedly solving instances of the original problem \eqref{eq:bMILP}, for each sample $i \in [n]$, points can be obtained and added to $\mathcal{Z}_{i}^{(0)}$. Specifically, in the present work, the first $n_{init}$ least costly solutions are included in $\mathcal{Z}_{i}^{(0)}$ for each $i \in [n]$. This is done by solving \eqref{eq:bMILP} $n_{init}$ times, using $\mathbf{y}_i$ as a cost vector, for each $i \in [n]$. Then, to get the $j+1$--\textit{th} least costly solution, no--good cuts of the form $\mathbf{z} \notin \mathcal{Z}_{i}^{(0)}$, are added to \eqref{eq:bMILP}, where $\mathcal{Z}_{i}^{(0)}$ already contains the first $j$ least costly solutions.
On the other hand, feasible solutions of \eqref{eq:bzV2} can be given as a warmstart to the solver, using any value of $\mathbf{w}$. This can be done by setting $\bar{\mathbf{z}}_i = \mathbf{z}_{i}^{SP} (\mathbf{w})$, and $\theta_i = \bar{\mathbf{y}}_{i}^{\intercal} \mathbf{z}_{i}^{SP}(\mathbf{w})$, for all $i \in [n]$, which results in a feasible point of \eqref{eq:bzV3}, whereas $\lambda_{ie}$ initial values are set accordingly. Notice that this setting permits using LR or SPO+ solutions (or any other) as a warmstart to calculate the IEO solution.


\subsubsection{Maximum number of points allowed} On the branch--and--cut scheme proposed in  Section~\ref{ssec:bac_gen}, the $\lambda_{ie}$ variables need to be defined for all $ i \in [n]$ and $e \in \mathcal{F}(\mathcal{Z})$. However, in general, $|\mathcal{F}(\mathcal{Z})|$ is a huge number, prohibiting a practical implementation. To that end, a lower number $n_{EP}$ of points can be used, defining $e \in [n_{EP}]$, and limiting the number of variables $\lambda_{ie}$. To validate then the correctness of the solution obtained, it is necessary to check whether the number of points $e(i) < n_{EP}$ at the end of the algorithm. If so, the provided solution is correct. If not,  $n_{EP}$ must be increased.

\subsubsection{Generation of multiple valid cuts}
Finally, multiple cuts can be added simultaneously by leveraging the information available at the cut generation stage of the branch-and-cut algorithm (Section~\ref{ssec:bac_gen}). Specifically, whenever a solution is obtained from the linear relaxation of a node, a value $\mathbf{w}^j$ is computed, from which a point in $\mathcal{Z}$ can be determined by solving \eqref{eq:spo_sp}. To enhance this process, we propose generating multiple points by perturbing $\mathbf{w}^j$. These perturbed vectors are defined as $\mathbf{w} = \mathbf{w}^j + \gamma \pmb{\epsilon}$, where $\gamma$ is a non-negative scalar, and $\pmb{\epsilon}$ is a uniformly distributed random vector in $[0, 1]^{n_x}$. As a result, up to $n_{rep}$ distinct values of $\mathbf{w}$ can be generated during each cut generation stage of the algorithm presented in Section~\ref{ssec:bac_gen}, potentially leading to the addition of up to $n_{rep}$ new cuts \eqref{eq:bac_cuts} to the formulation in a single iteration.

\section{Computational experiments}

In this section, we conduct computational experiments to evaluate the performance of our pessimistic IEO loss minimization algorithm for combinatorial optimization problems. Specifically, we compare the $\ell_{IEO}$ loss with the mean squared error ($\ell_{MSE}$) and the SPO+ loss ($\ell_{SPO+}$) proposed in \cite{Elmachtoub2022}. The experiments involve multiple instances of the 0-1 knapsack problem, varying both the randomness generation process and the number of items.

For clarity, the loss functions under consideration are defined as follows:
\begin{subequations} \label{eq:exp_losses}
\begin{align}
    \ell_{MSE}(f(\mathbf{x}), \mathbf{y}) &= || f (\mathbf{x}_i) - \mathbf{y}_i ||^2, \label{eq:exp_losses_lr} \\
    \ell_{SPO+}(f(\mathbf{x}), \mathbf{y}) &= \underset{\mathbf{z} \in {\mathcal{Z}}}{\text{max }} \lbrace \mathbf{y}^{\intercal} \mathbf{z} - 2 f(\mathbf{x})^{\intercal} \mathbf{z} \rbrace + 2 f(\mathbf{x})^{\intercal} \mathbf{z}^{\star} - \mathbf{y}^{\intercal} \mathbf{z}^{\star}. \label{eq:exp_losses_spop}
\end{align}
\end{subequations}

The $\ell_{MSE}$ loss does not incorporate any information about the DOP, whereas the SPO+ loss ($\ell_{SPO+}$) is an IEO method that serves as an approximation of the $\ell_{IEO}$ loss. Hereinafter, we refer to the calibration methods used in ERM as follows: LR for $\ell_{MSE}$, SPO+ for $\ell_{SPO+}$, and IEO for $\ell_{IEO}$.


We perform out-of-sample experiments in which different predictors are calculated from the same datasets by minimizing the three previously mentioned loss functions. We consider linear predictors of the form $f_{\mathbf{w}}(\mathbf{x}) = \mathbf{W} \mathbf{x}$, which along with the linearization \eqref{eq:mc_relax}, renders the model a MILP one, allowing the use of off--the--shelf solvers. Regarding the ERM resolution, for $\ell_{MSE}$ we use available packages, whereas in the case of $\ell_{SPO+}$, we implemented the solution approach presented in \cite{Elmachtoub2022}, Section 5.1. It is worth mentioning that the solution methods presented in \cite{Elmachtoub2022} for solving $\ell_{SPO+}$ minimization are for continuous optimization problems. The authors suggest to use the convex hull description for integer optimization problems. Unfortunately, this is rarely practical due to the large number of inequalities required to describe the convex hull.
As a result, we use the continuous relaxation of the knapsack problem in minimizing the SPO+ loss, 
but, to see the impact of using this relaxation, we also implemented SPO+ using the convex hull description for small size instances. In the case of the ERM of $\ell_{IEO}$, the branch--and--cut algorithm, along with the settings presented in Section~\ref{ssec:bac_gen} and  Section~\ref{ssec:settings}, respectively, is used. Finally, out--of--sample experiments are performed over a dataset of unseen samples $\mathcal{D}_{test} = \lbrace (\mathbf{x}_i , \mathbf{y}_{i}) \rbrace_{i=1}^{n_{test}}$.

All optimization models were solved using Gurobi 11.0.0. Unless otherwise specified, a time limit of 1800 seconds was imposed on MILP models. Computational experiments were conducted on a system equipped with an Intel Core i7-10870H processor (2.2 GHz) and 16 GB of RAM.

\subsection{0--1 Knapsack problem instances:}
For the experiments, we create instances of the 0-1 knapsack problem using the same methodology as the one presented in \cite{Fatma2022}, Section 6.2.
The weights of the items, $a_j$ for $j \in [n_y]$, are random integer values between 1 and 1000, whereas the knapsack capacity $b$ is a random integer between $b_{l}$ and $b_{u}$, where $b_{l} = \text{max}_{j \in [n_y]} a_j$, and $b_{u} = (1 - (1 - r)\frac{b_l}{\mathbf{1}^{\intercal} \mathbf{a}})\mathbf{1}^{\intercal} \mathbf{a}$, where $r$ is uniformly distributed on $[0, 1]$. The parameters to be estimated are the objective function coefficients of the items, and they are associated with five features ($n_x = 5$). We have a ground truth matrix $W_{GT} \in \mathbb{R}^{n_y \times n_x}$ associated with each instance, with coefficients taken randomly. Each component of the feature vector $\mathbf{x}_i$ is taken from a uniform distribution on $[-1, 1]$, except for the last entry which is set equal to $1$, to model a constant term on the predictor. In that way, each component of the cost vector $\mathbf{y}_i$ is generated as:
\begin{equation*}
    y_{ij} = \Tilde{\alpha}_{ij} \big ( \mathbf{w}_{GT, j}^{\intercal} \mathbf{x}_{i} \big )^{\delta} + \eta_{ij} \quad j \in [n_y], \, i \in [n],
\end{equation*}
where $\mathbf{w}_{GT, j}$ is the $j$--th row of $W_{GT}$, $\Tilde{\alpha}_{ij}$ is uniformly distributed on $[1-\alpha, 1+\alpha]$, and $2 \eta_{ij} + 1$ is an exponential random variable with scale parameter $\lambda = 1$.
Following \cite{Fatma2022}, a normalized version of the IEO loss is used as a measure of out-of-sample performance:
\begin{equation*}
    \bar{\ell}_{IEO} (f(\mathbf{x}), \mathbf{y}) = \dfrac{\underset{\mathbf{z} \in \mathcal{Z}}{\text{max }} c (\mathbf{y}, \mathbf{z} )  -  \underset{\mathbf{z} \in \mathcal{Z}^\star (f_{\mathbf{w}}(\mathbf{x}))}{\text{min }} c (\mathbf{y}, \mathbf{z} )}{\underset{\mathbf{z} \in \mathcal{Z}}{\text{max }} c (\mathbf{y}, \mathbf{z} )},
\end{equation*}
where $\% \: \bar{\ell}_{IEO}$ will denote the percentage value of $\bar{\ell}_{IEO}$.

\subsection{Comparison between CCG and branch--and--cut algorithms}
We start by testing the performance of the proposed branch-and-cut algorithm against that of the CCG algorithm. We do so by measuring the in-sample $\bar{\ell}_{IEO}$  after a defined time limit of computation, over two of the instances used in \cite{Fatma2022} for the knapsack problem.
In the case of the CCG algorithm, it is also required to set a time limit to the master problem resolution $P_{\varepsilon}^{(k)}$. We set this to  $1/10$ of the algorithm's time limit.

In Table \ref{tab:times_alg}, the final in-sample $\% \: \bar{\ell}_{IEO}$ value is shown for both algorithms, as well as the number of pairs of cuts \eqref{eq:bzV3_2}--\eqref{eq:bzV3_4} added (\# cuts) and the final gap (gap) over a set of different instances. These are constructed by setting time limits of 600, 1200, and 1800s over two instances of the knapsack problem, using $n \in \{100, 300, 500 \}$. The first observation is that the branch-and-cut algorithm allows us to get predictors with lower in-sample $\bar{\ell}_{IEO}$ (lower values are in boldface).
Additionally, in two instances (Instance 10, with $n = 100$ and $n = 500$), the branch-and-cut algorithm solves the problem to optimality, achieving a gap smaller than the predefined tolerance ($10^{-4}$). In contrast, the CCG algorithm does not reach this level of precision.

\begin{table}[t!]
\centering
\begin{tabular}{ccc|ccc|ccc}
         &     &            & \multicolumn{3}{c|}{CCG}           & \multicolumn{3}{c}{B\&C}           \\ \noalign{\smallskip}\hline\noalign{\smallskip}
Instance & $n$ & Time limit & \# cuts & $\% \: \bar{\ell}_{IEO}$ & gap & \# cuts & $\% \: \bar{\ell}_{IEO}$ & gap \\ \noalign{\smallskip}\hline\noalign{\smallskip}
1 & 100 & 600s & 5  & 0.193 & 0.187 & 656 & \textbf{0.027} & 0.021 \\
1 & 300 & 600s & 5  & 0.318 & 0.318 & 216 & \textbf{0.182} & 0.312 \\
1 & 500 & 600s & 4  & 0.358 & 0.357 & 774 & \textbf{0.358} & 0.347 \\
10 & 100 & 600s & 10 & 0.015 & 0.015 & 660 & \textbf{0.007} & $<10^{-4}(\star)$ \\
10 & 300 & 600s & 0  & 0.203 & 0.202 & 582 & \textbf{0.041} & 0.042 \\
10 & 500 & 600s & 1  & 0.227 & 0.225 & 3   & \textbf{0.031} & $<10^{-4}(\star)$ \\
1 & 100 & 1200s & 3  & 0.223 & 0.171 & 818 & \textbf{0.027} & 0.021 \\
1 & 300 & 1200s & 5  & 0.318 & 0.318 & 1228 & \textbf{0.181} & 0.312 \\
1 & 500 & 1200s & 4  & 0.358 & 0.355 & 784 & \textbf{0.224} & 0.347 \\
10 & 100 & 1200s & 10 & 0.015 & 0.015 & 660 & \textbf{0.007} & $<10^{-4}(\star)$ \\
10 & 300 & 1200s & 0  & 0.203 & 0.202 & 680 & \textbf{0.032} & 0.034 \\
10 & 500 & 1200s & 1  & 0.227 & 0.225 & 3   & \textbf{0.031} & $<10^{-4}(\star)$ \\
1 & 100 & 1800s & 7  & 0.205 & 0.132 & 944 & \textbf{0.026} & 0.020 \\
1 & 300 & 1800s & 5  & 0.318 & 0.318 & 1764 & \textbf{0.148} & 0.312 \\
1 & 500 & 1800s & 4  & 0.358 & 0.355 & 804 & \textbf{0.185} & 0.347 \\
10 & 100 & 1800s & 10 & 0.015 & 0.015 & 660 & \textbf{0.007} & $<10^{-4}(\star)$ \\
10 & 300 & 1800s & 0  & 0.203 & 0.202 & 688 & \textbf{0.032} & 0.034 \\
10 & 500 & 1800s & 1  & 0.227 & 0.225 & 3   & \textbf{0.031} & $<10^{-4}(\star)$ \\ \noalign{\smallskip}\hline
\end{tabular}
\caption{Computational performance comparison between CCG and branch-and-cut algorithms. 0--1 Knapsack instances with $n_y=10$, $\delta = 3$. $(\star)$ Solved before reaching the time limit.}
\label{tab:times_alg}
\end{table}

One of the most salient features of Table \ref{tab:times_alg} is the significant difference in the number of added cuts between the algorithms, which can be attributed as the main cause behind the differences in the $\bar{\ell}_{IEO}$ values. This difference is explained by the fact that the CCG algorithm needs to solve the master problem to generate new cuts, which can be very time consuming. Besides, the time limit imposed on the master problem could be insufficient to allow it to compute new cuts. Although early-stopping settings can be adjusted to expedite the CCG algorithm, its performance remains highly instance-dependent. For instance, in the case of Instance 1 with $n=100$, the CCG algorithm fails to achieve lower in-sample $\bar{\ell}_{IEO}$ values despite an increased time limit. This outcome occurs because the master problem reaches its time limit in these scenarios, generating points that are already included in the sets $\mathcal{Z}_{i}^{(k)}$, thereby reducing the number of cuts added. Consequently, if none of the sets $\mathcal{Z}_{i}^{(k)}$ are updated, an early-stopping condition is triggered, halting the algorithm.

In contrast, configurations for the branch-and-cut algorithm are simpler, considering that built-in functions are available in many off-the-shelf solvers, allowing us to efficiently add cuts \eqref{eq:bzV3_2}--\eqref{eq:bzV3_4} as lazy constraints.

\subsection{In--sample v/s out--of--sample performance and overfitting} \label{ssec:overfitting}

In this section, we evaluate the quality of the in-sample loss as an estimator of the out-of-sample loss. The performance metric used to evaluate predictors is the normalized pessimistic IEO loss, $\bar{\ell}_{IEO}$, which is widely employed in the context of CO. To the best of the authors' knowledge, the proposed method is the first to directly minimize the $\bar{\ell}_{IEO}$ loss function to calibrate predictors. This approach is expected to yield better in-sample performance compared to other methods. For example, the SPO+ loss \eqref{eq:exp_losses_spop} serves as a surrogate for $\ell_{IEO}$, and significant approximations are employed in \cite{Bucarey2024} to minimize $\ell_{IEO}$.

Table \ref{tab:init_ov1} presents the in-sample and out-of-sample $\% \: \bar{\ell}_{IEO}$ results for the knapsack problem with $n_y = 10$ and $n \in \{100, 200\}$, comparing three studied methods. Additionally, we include SPO+(CH), a variant of SPO+ that incorporates the convex hull, computed using the PORTA software \cite{porta_1.4.1}. It is worth noting that SPO+(CH) is feasible only for low-dimensional instances, where the convex hull description is obtainable.

From Table \ref{tab:init_ov1}, we observe that IEO consistently outperforms LR, SPO+, and SPO+(CH) in in-sample performance. This is expected, as IEO is the only method among the four that directly minimizes the in-sample $\bar{\ell}_{IEO}$. However, in the out-of-sample experiments, no method consistently outperforms the others.

\begin{table}[t!]
\centering
\resizebox{\textwidth}{!}{
\begin{tabular}{ccc|cccc|cccc}
\multicolumn{1}{l}{}  &          &       & \multicolumn{4}{c|}{In--sample}                 & \multicolumn{4}{c}{Out--of--sample}                               \\ \noalign{\smallskip}\hline\noalign{\smallskip}
\multicolumn{1}{l}{$n$} & Instance & $\delta$ & LR    & SPO+  & SPO+(CH)       & IEO            & LR             & SPO+           & SPO+(CH)       & IEO            \\ \noalign{\smallskip}\hline\noalign{\smallskip}
100                  & 0        & 1     & 1.023 & 1.231 & 0.892          & \textbf{0.820} & 0.954          & 1.093          & \textbf{0.951} & 1.056          \\
100                  & 0        & 3     & 0.540 & 0.298 & 0.029          & \textbf{0.015} & 0.458          & 0.441          & \textbf{0.429} & 0.486          \\
100                  & 0        & 5     & 0.993 & 0.325 & \textbf{0.000} & \textbf{0.000} & 1.020          & 0.828          & \textbf{0.680} & \textbf{0.680} \\
100                  & 5        & 1     & 3.583 & 3.244 & 3.320          & \textbf{3.197} & \textbf{3.013} & 3.360          & 3.144          & 3.302          \\
100                  & 5        & 3     & 0.227 & 0.209 & \textbf{0.017} & \textbf{0.017} & \textbf{0.533} & 0.591          & 0.819          & 0.819          \\
100                  & 5        & 5     & 1.380 & 0.468 & \textbf{0.006} & \textbf{0.006} & 1.741          & \textbf{1.044} & 1.193          & 1.193          \\
100                  & 7        & 1     & 0.475 & 0.341 & 0.413          & \textbf{0.222} & 0.726          & \textbf{0.721} & 0.724          & 1.960          \\
100                  & 7        & 3     & 0.477 & 0.149 & 0.049          & \textbf{0.024} & 0.318          & 0.286          & \textbf{0.212} & 0.706          \\
100                  & 7        & 5     & 1.445 & 0.534 & \textbf{0.004} & \textbf{0.004} & 1.483          & 1.053          & \textbf{0.854} & \textbf{0.854} \\
200                  & 0        & 1     & 1.025 & 1.119 & 0.893          & \textbf{0.863} & \textbf{0.996} & 1.140          & 1.036          & 1.101          \\
200                  & 0        & 3     & 0.518 & 0.270 & 0.044          & \textbf{0.024} & 0.420          & 0.454          & \textbf{0.222} & 0.300          \\
200                  & 0        & 5     & 0.927 & 0.409 & \textbf{0.044} & \textbf{0.044} & 0.923          & 0.558          & \textbf{0.262} & \textbf{0.262} \\
200                  & 5        & 1     & 3.685 & 3.442 & 3.592          & \textbf{3.279} & \textbf{2.982} & 3.473          & 3.119          & 3.407          \\
200                  & 5        & 3     & 0.483 & 0.410 & \textbf{0.255} & \textbf{0.255} & 0.549          & 0.664          & \textbf{0.425} & \textbf{0.425} \\
200                  & 5        & 5     & 1.225 & 0.936 & \textbf{0.154} & \textbf{0.154} & 1.366          & 0.887          & \textbf{0.698} & \textbf{0.698} \\
200                  & 7        & 1     & 0.707 & 0.690 & 0.651          & \textbf{0.642} & 0.813          & 0.824          & \textbf{0.797} & 0.804          \\
200                  & 7        & 3     & 0.457 & 0.318 & 0.091          & \textbf{0.070} & 0.323          & \textbf{0.217} & 0.368          & 0.459          \\
200                  & 7        & 5     & 1.419 & 0.729 & \textbf{0.098} & \textbf{0.098} & 1.291          & 1.076          & \textbf{0.701} & \textbf{0.701} \\ \noalign{\smallskip}\hline
\end{tabular}
}
\caption{In--sample v/s Out--of--sample $\% \: \bar{\ell}_{IEO}$ comparison. 0--1 Knapsack instances with $n_y=10$ and $\alpha = 0.1$, tested over $n_{test} = 1000$ unseen samples.}
\label{tab:init_ov1}
\end{table}

\begin{table}[t!]
\centering
\resizebox{\textwidth}{!}{
\begin{tabular}{ccc|cccc|cccc}
\multicolumn{1}{l}{}  &          &       & \multicolumn{4}{c|}{In--sample}                 & \multicolumn{4}{c}{Out--of--sample}                               \\ \noalign{\smallskip}\hline\noalign{\smallskip}
\multicolumn{1}{l}{$n$} & Instance & $\delta$ & LR    & SPO+  & SPO+(CH)       & IEO            & LR             & SPO+           & SPO+(CH)       & IEO            \\ \noalign{\smallskip}\hline\noalign{\smallskip}
100                  & 0        & 1     & 0.736          & 1.357 & 0.826          & \textbf{0.721} & \textbf{1.261} & 1.383 & 1.349          & 1.334          \\
100                  & 0        & 3     & 0.608          & 0.443 & 0.024          & \textbf{0.021} & 0.809          & 0.653 & 0.371          & \textbf{0.369} \\
100                  & 0        & 5     & 1.280          & 0.301 & \textbf{0.000} & \textbf{0.000} & 1.620          & 0.973 & \textbf{0.870} & \textbf{0.870} \\
100                  & 5        & 1     & 2.892          & 3.031 & 2.892          & \textbf{1.909} & \textbf{3.721} & 4.321 & 3.942          & 4.684          \\
100                  & 5        & 3     & 0.315          & 0.400 & \textbf{0.012} & \textbf{0.012} & \textbf{0.794} & 1.013 & 1.082          & 1.082          \\
100                  & 5        & 5     & 1.242          & 1.130 & 0.147          & \textbf{0.085} & 2.247          & 1.974 & \textbf{1.246} & 1.274          \\
100                  & 7        & 1     & 0.578          & 0.724 & 0.534          & \textbf{0.308} & 0.948          & 1.036 & 0.981          & \textbf{0.945} \\
100                  & 7        & 3     & 0.544          & 0.097 & 0.041          & \textbf{0.024} & 0.598          & 0.425 & 0.316          & \textbf{0.310} \\
100                  & 7        & 5     & 2.222          & 0.550 & \textbf{0.000} & \textbf{0.000} & 1.780          & 1.662 & \textbf{1.097} & \textbf{1.097} \\
200                  & 0        & 1     & 0.816          & 0.901 & 0.830          & \textbf{0.791} & 1.226          & 1.368 & \textbf{1.188} & 1.317          \\
200                  & 0        & 3     & 0.474          & 0.496 & 0.029          & \textbf{0.014} & 0.667          & 0.628 & \textbf{0.195} & 0.200          \\
200                  & 0        & 5     & 1.205          & 0.584 & 0.070          & \textbf{0.056} & 1.538          & 0.916 & 0.505          & \textbf{0.487} \\
200                  & 5        & 1     & \textbf{3.012} & 3.317 & 3.393          & \textbf{3.012} & \textbf{3.438} & 3.746 & 3.572          & \textbf{3.438} \\
200                  & 5        & 3     & 0.641          & 0.663 & 0.108          & \textbf{0.081} & 0.694          & 0.967 & \textbf{0.555} & \textbf{0.555} \\
200                  & 5        & 5     & 1.140          & 0.744 & 0.081          & \textbf{0.026} & 1.775          & 1.650 & 0.877          & \textbf{0.868} \\
200                  & 7        & 1     & 0.609          & 0.531 & 0.486          & \textbf{0.388} & \textbf{0.938} & 1.008 & 0.984          & 1.025          \\
200                  & 7        & 3     & 0.335          & 0.278 & 0.086          & \textbf{0.061} & \textbf{0.454} & 0.510 & 0.802          & 0.810          \\
200                  & 7        & 5     & 1.295          & 0.989 & 0.046          & \textbf{0.040} & 1.743          & 2.032 & \textbf{0.504} & \textbf{0.504} \\ \noalign{\smallskip}\hline
\end{tabular}
}
\caption{In--sample v/s Out--of--sample $\% \: \bar{\ell}_{IEO}$ comparison using LHS and lasso regularization. 0--1 Knapsack instances with $n_y=10$ and $\alpha = 0.1$, tested over $n_{test} = 1000$ unseen samples.}
\label{tab:init_ov2}
\end{table}

Usually, the in-sample loss value is expected to be a decent proxy for the out-of-sample loss. However, this is not always the case, and when in-sample performance significantly underestimates the out-of-sample behavior, we are in the presence of what is commonly known as overfitting. In the context of prediction, overfitting is a well-studied phenomenon, see, e.g., \cite{bejani2021systematic} and \cite{montesinos2022overfitting}.
In contrast, in IEO problems overfitting is not a well-understood subject. In a recent work, \cite{bennouna2023} studied it in the context of data-driven optimization, identifying three sources of overfitting: statistical error, the presence of noise, and also misspecification or the presence of contaminated data. The first one refers to the effect of having limited data samples, which can be insufficient to accurately capture the true underlying distribution, whereas the last one is related with the presence of corrupted data. Therefore, in  our case only the statistical error and the noise effect are plausible causes.

The study of consistent techniques to mitigate overfitting for predictive-prescriptive estimators is out of the scope of this work, and probably deserves a separate analysis. However, we test two remedies to cope with it. Firstly, under the hypothesis that the misidentification of the underlying distribution has a bigger impact for IEO methods (SPO+, SPO+(CH), IEO), we adopt latin hypercube sampling (LHS). LHS is a type of stratified sampling, in which the domain of the features parameters $\mathcal{X}$ is subdivided into $n_{div}$ non-overlapping subsets $\lbrace \mathcal{X}_p \rbrace_{p=1}^{n_{div}}$, where samples are taken from every subset with equal probability. This technique allows to uniformly cover the entire domain $\mathcal{X}$, helping to cope with the statistical error. Secondly, lasso regularization is applied to IEO, i.e., the term $\gamma \sum_{j=1}^{n_y} \sum_{l=1}^{n_x} |W_{jl}|$ is added to the objective function, with $\gamma = 0.005$.

Table \ref{tab:init_ov2} shows the results after applying LHS to all four methods, and lasso regularization to IEO, for the same experiments as the ones present in Table \ref{tab:init_ov1}. Firstly, IEO keeps outperforming the rest of the methods in-sample. We observe significant improvements in the out-of-sample performance of IEO with respect to the cases of Table \ref{tab:init_ov1}, with the lowest out-of-sample $\bar{\ell}_{IEO}$ in 10 of the 18 instances. Nevertheless, some pathological cases are still present that deserve further study. In Instance 5 with $n=100$ and $\delta=1$, the in-sample $\bar{\ell}_{IEO}$ values improve for all of the four methods, while the out-of-sample loss gets worse for all of them, with IEO being the best in-sample and the worst out-of-sample method. Then, taking the same case but considering $n=200$, we note that the differences are considerably reduced.

\subsection{Out-of-sample experiments}
We now present extensive computational experiments on the 0-1 knapsack problem.
In the first series of experiments, eleven instances of $n_y = 10$ are tested, taken directly from the ones published in \cite{Fatma2022} without modifications. In this case, SPO+(CH) is also included since the size of the instances allows the computation of the convex hull. As before, we consider $\delta \in \lbrace 1, 3, 5 \rbrace$, $\alpha = 0.1$, whereas the number of training  samples is selected as $n \in \lbrace 100, 200, 300, 400, \allowbreak 500  \rbrace$. Out-of-sample tests are made over $n_{test} = 1000$ unobserved samples, reporting the mean $\% \: \bar{\ell}_{IEO}$ value.

\begin{figure}
    \centering
    \includegraphics[width=1\linewidth]{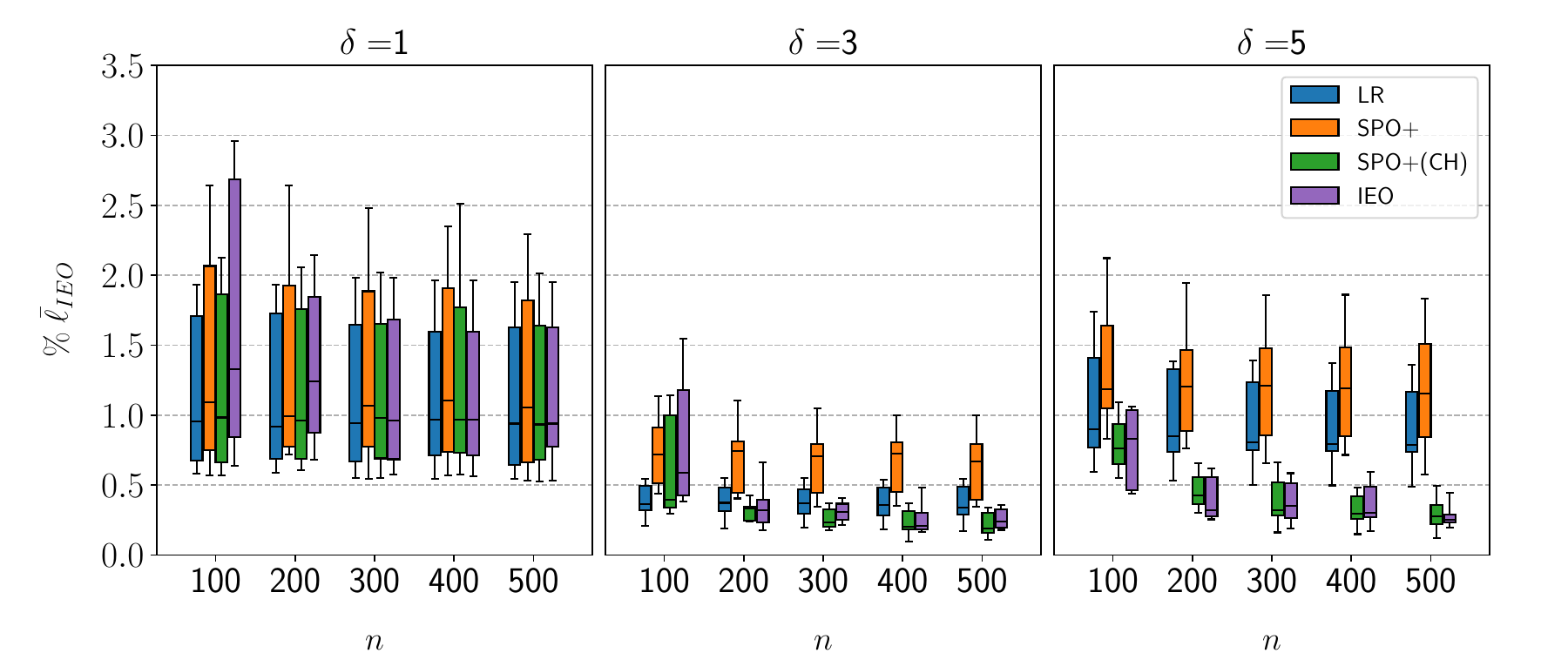}
    \caption{Out--of--sample results, $n_y = 10$, $n_{test} = 1000$.}
    \label{fig:ni10_nt_1000_lhs}
\end{figure}

In Figure~\ref{fig:ni10_nt_1000_lhs}, boxplots of the out-of-sample mean $\% \: \bar{\ell}_{IEO}$ values are shown for the four benchmarked methods. For $\delta = 1$, the performances are similar for the four methods, with LR being slightly better that its peers. As we move away from linearity, with $\delta = 3$ or $5$, a clear dominance of SPO+(CH) and IEO over LR and SPO+ is observed. We highlight that SPO+'s performance is the worst in those cases, which can be attributed to the quality of the linear relaxation of the knapsack problem used by that method. Additionally, along the same lines of the observations made in Section \ref{ssec:overfitting}, for a lower number of samples ($n = 100$), IEO overfits for $\delta = 1$. However, as the number of training samples $n$ increases, it performs practically the same as LR and SPO+(CH).

To evaluate the performance of different methods as the size of the knapsack problem increases, we also ran experiments on instances with $n_y = 20$. Specifically, 10 instances are generated for each value of $n_y$ using the methodology described in \cite{Fatma2022}, with $\delta \in \{1, 3, 5\}$ applied to each instance. For these experiments, SPO+(CH) could not be tested due to the inability to obtain the convex hull description of $\mathcal{Z}$.

\begin{figure}
  \includegraphics[width=1\linewidth]{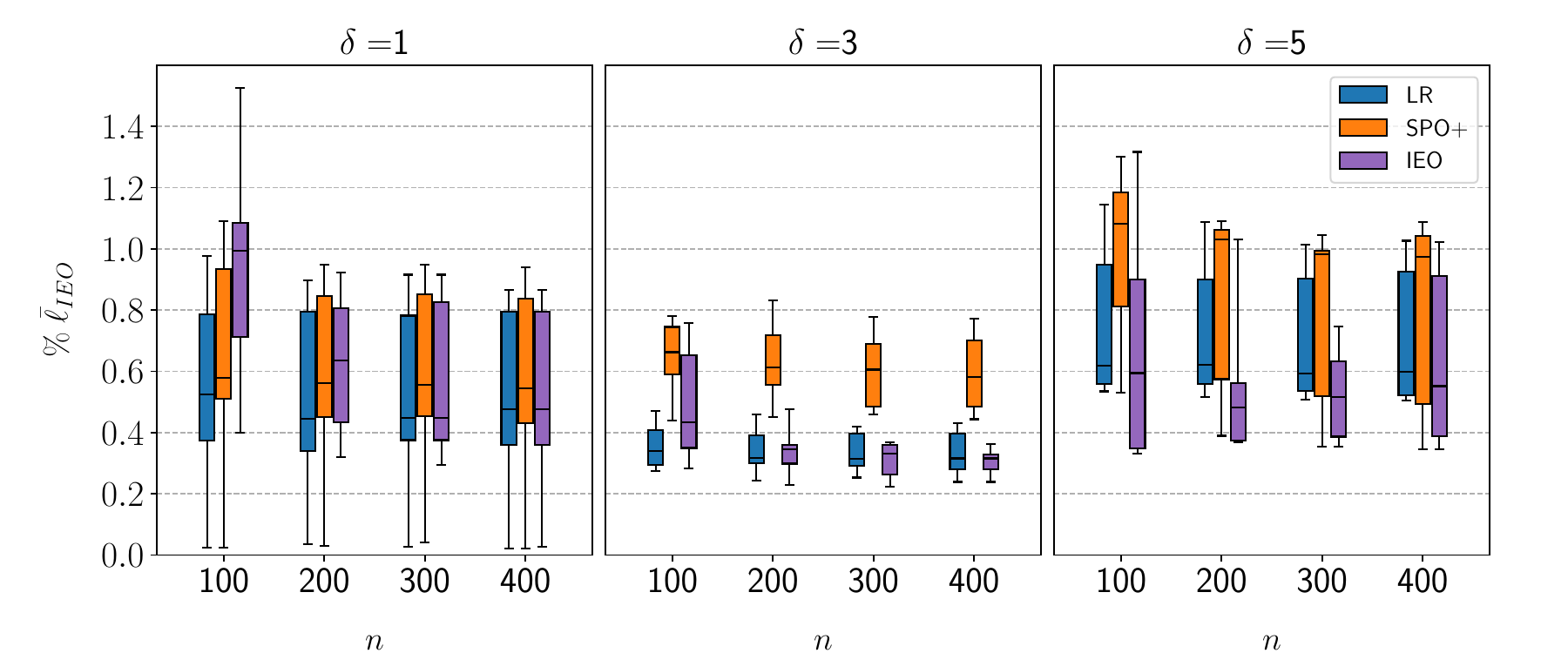}
\caption{Out--of--sample results, $n_y = 20$, $n_{test} = 1000$.}
\label{fig:ni20_nt_1000_lhs}  
\end{figure}

Figure~\ref{fig:ni20_nt_1000_lhs} presents boxplots of the out-of-sample mean values $\% \: \bar{\ell}_{IEO}$, computed using 1000 unobserved samples. The number of training samples was restricted to $n \in \{100, 200, 300, 400\}$ (the case with $n=500$ was excluded), as these were the cases where IEO successfully improved the warm-start solutions provided to the solver within a time limit of 1800 seconds. Overall, the results exhibit a behavior consistent with that observed in Figure~\ref{fig:ni10_nt_1000_lhs}. For $n = 100$, overfitting is evident, particularly for $\delta = 1$. However, as the number of training samples increases, this overfitting effect diminishes, and the three methods demonstrate similar performance. For $\delta = 3$, as with the case of $n_y = 10$, SPO+ shows the worst performance, while LR and IEO achieve comparable average $\bar{\ell}_{IEO}$ loss values. Finally, for $\delta = 5$, IEO outperforms both LR and SPO+, with SPO+ showing the poorest performance.
It is worth noting that, for $n = 400$, similar results were observed due to IEO's difficulty in improving the warmstart solutions within the given time limit, particularly in terms of the objective function value. This limitation suggests that the algorithm may require more time to refine the initial solutions, indicating a potential area for further research to enhance its computational efficiency. Nevertheless, across the analyzed cases, IEO demonstrates robustness against estimator function misspecification and strong performance when the predictor is well specified (cases where $\delta = 1$).

\section{Conclusions and future research}
In this work, we propose a novel method for producing IEO estimators by minimizing the normalized $\ell_{IEO}$ loss function, avoiding the use of surrogate functions or approximations. Unlike previous approaches, our method enables the treatment of downstream optimization problems
containing binary variables without requiring the convex hull description of the feasible sets. To solve this problem, we introduce two algorithms: a CCG algorithm, inspired by \cite{BoZeng2020}, and a branch-and-cut algorithm.

We conducted computational experiments to evaluate both the computational and out-of-sample performance of the proposed approach. Our findings show that the branch-and-cut algorithm outperforms the CCG algorithm in terms of the final in-sample loss (objective function value), enabling the computation of estimators with lower in-sample $\bar{\ell}_{IEO}$ values within the same time limit. Furthermore, numerical experiments on the 0–1 knapsack problem demonstrate that the proposed IEO method outperforms both LR and SPO+ in terms of out-of-sample IEO
loss. Additionally, IEO exhibits robustness against model misspecification, unlike LR and SPO+.

When the convex hull description is provided for its formulations, SPO+ behaves similarly to IEO. However, the requirement for a convex hull description is a significant limitation when dealing with DOPs containing a moderate number of integer variables.
As a result, SPO+ exhibited poor out-of-sample performance when using the fractional knapsack (linear relaxation)-based estimator, performing even worse than the LR-based estimator, which completely disregards DOP information.

Finally, we observed that when using a small number of samples, the proposed IEO method, as well as SPO+ with the convex hull description, suffers from overfitting, even when achieving the lowest in-sample $\bar{\ell}_{IEO}$ loss values. We tested two simple measures to mitigate this issue, with some success in reducing its effects, though not completely eliminating them. This suggests a potential pathological behavior in IEO methods that warrants further investigation.

\bibliographystyle{ieeetr}    
\bibliography{bibtex.bib}   

\begin{thebibliography}{10}

\bibitem{LIU2023}
Z.~Liu, Y.~Yin, F.~Bai, and D.~K. Grimm, ``End-to-end learning of user
  equilibrium with implicit neural networks,'' {\em Transportation Research
  Part C: Emerging Technologies}, vol.~150, p.~104085, 2023.

\bibitem{islam2023}
M.~S. Islam and A.~T. Wasi, ``Optimizing inventory routing: A decision-focused
  learning approach using neural networks,'' {\em arXiv preprint
  arXiv:2311.00983}, 2023.

\bibitem{xu-cohen-2018}
Y.~Xu and S.~B. Cohen, ``Stock movement prediction from tweets and historical
  prices,'' in {\em Proceedings of the 56th Annual Meeting of the Association
  for Computational Linguistics (Volume 1: Long Papers)} (I.~Gurevych and
  Y.~Miyao, eds.), (Melbourne, Australia), pp.~1970--1979, Association for
  Computational Linguistics, 2018.

\bibitem{pagnoncelli2023synthetic}
B.~K. Pagnoncelli, D.~Ram{\'\i}rez, H.~Rahimian, and A.~Cifuentes, ``A
  synthetic data-plus-features driven approach for portfolio optimization,''
  {\em Computational Economics}, vol.~62, no.~1, pp.~187--204, 2023.

\bibitem{wang2023}
S.~Wang, Y.~Bai, T.~Ji, K.~Fu, L.~Wang, and C.-T. Lu, ``Stock movement and
  volatility prediction from tweets, macroeconomic factors and historical
  prices,'' {\em arXiv preprint arXiv:2312.03758}, 2023.

\bibitem{pena2024modified}
J.-M. Pe{\~n}a, F.~Su{\'a}rez, O.~Larr{\'e}, D.~Ram{\'\i}rez, and A.~Cifuentes,
  ``A modified ctgan-plus-features-based method for optimal asset allocation,''
  {\em Quantitative Finance}, vol.~24, no.~3-4, pp.~465--479, 2024.

\bibitem{sadana2023survey}
U.~Sadana, A.~Chenreddy, E.~Delage, A.~Forel, E.~Frejinger, and T.~Vidal, ``A
  survey of contextual optimization methods for decision making under
  uncertainty,'' {\em arXiv preprint arXiv:2306.10374}, 2023.

\bibitem{Bertsimas2020}
D.~Bertsimas and N.~Kallus, ``From predictive to prescriptive analytics,'' {\em
  Management Science}, vol.~66, no.~3, pp.~1025--1044, 2020.

\bibitem{elmachtoub2023estimate}
A.~N. Elmachtoub, H.~Lam, H.~Zhang, and Y.~Zhao, ``Estimate-then-optimize
  versus integrated-estimation-optimization versus sample average
  approximation: A stochastic dominance perspective,'' {\em arXiv preprint
  arXiv:2304.06833}, 2023.

\bibitem{mandi2024}
J.~Mandi, J.~Kotary, S.~Berden, M.~Mulamba, V.~Bucarey, T.~Guns, and
  F.~Fioretto, ``Decision-focused learning: Foundations, state of the art,
  benchmark and future opportunities,'' 2024.

\bibitem{Malaga2022}
M.~Muñoz, S.~Pineda, and J.~Morales, ``A bilevel framework for decision-making
  under uncertainty with contextual information,'' {\em Omega}, vol.~108,
  p.~102575, 2022.

\bibitem{ban2019big}
G.-Y. Ban and C.~Rudin, ``The big data newsvendor: Practical insights from
  machine learning,'' {\em Operations Research}, vol.~67, no.~1, pp.~90--108,
  2019.

\bibitem{Elmachtoub2022}
A.~N. Elmachtoub and P.~Grigas, ``Smart “predict, then optimize”,'' {\em
  Management Science}, vol.~68, no.~1, pp.~9--26, 2022.

\bibitem{Bucarey2024}
V.~Bucarey, S.~Calder{\'o}n, G.~Mu{\~{n}}oz, and F.~Semet, ``Decision-focused
  predictions via pessimistic bilevel optimization: A computational study,''
  in {\em Integration of Constraint Programming, Artificial Intelligence, and
  Operations Research} (B.~Dilkina, ed.), (Cham), pp.~127--135, Springer Nature
  Switzerland, 2024.

\bibitem{Fatma2022}
N.~Ho-Nguyen and F.~K\i{}l\i{}n\c{c}-Karzan, ``Risk guarantees for end-to-end
  prediction and optimization processes,'' {\em Management Science}, vol.~68,
  no.~12, pp.~8680--8698, 2022.

\bibitem{Estes2023}
A.~S. Estes and J.-P.~P. Richard, ``Smart predict-then-optimize for two-stage
  linear programs with side information,'' {\em INFORMS Journal on
  Optimization}, 2023.

\bibitem{mandi2019smart}
J.~Mandi, E.~Demirović, P.~J. Stuckey, and T.~Guns, ``Smart
  predict-and-optimize for hard combinatorial optimization problems,'' {\em
  arXiv preprint arXiv:1911.10092}, 2019.

\bibitem{Wiesemann2013}
W.~Wiesemann, A.~Tsoukalas, P.-M. Kleniati, and B.~Rustem, ``Pessimistic
  bilevel optimization,'' {\em SIAM Journal on Optimization}, vol.~23, no.~1,
  pp.~353--380, 2013.

\bibitem{Goerigk2024}
M.~Goerigk, J.~Kurtz, M.~Schmidt, and J.~Thürauf, ``Connections and
  reformulations between robust and bilevel optimization,'' {\em Optmization
  online}, 2024.

\bibitem{Dempe2018}
S.~Dempe, G.~Luo, and S.~Franke, ``Pessimistic bilevel linear optimization,''
  {\em Journal of Nepal Mathematical Society}, vol.~1, p.~1–10, Feb. 2018.

\bibitem{KLEINERT2021}
T.~Kleinert, M.~Labbé, I.~Ljubić, and M.~Schmidt, ``A survey on mixed-integer
  programming techniques in bilevel optimization,'' {\em EURO Journal on
  Computational Optimization}, vol.~9, p.~100007, 2021.

\bibitem{Koppe2010}
M.~Köppe, M.~Queyranne, and C.~T. Ryan, ``Parametric integer programming
  algorithm for bilevel mixed integer programs,'' {\em Journal of Optimization
  Theory and Applications}, vol.~146, 2010.

\bibitem{Fischetti2017}
M.~Fischetti, I.~Ljubi\'{c}, M.~Monaci, and M.~Sinnl, ``A new general-purpose
  algorithm for mixed-integer bilevel linear programs,'' {\em Operations
  Research}, vol.~65, no.~6, pp.~1615--1637, 2017.

\bibitem{Tahernejad2020}
S.~Tahernejad, T.~K. Ralphs, and S.~T. DeNegre, ``A branch-and-cut algorithm
  for mixed integer bilevel linear optimization problems and its
  implementation,'' {\em Mathematical Programming Computation}, vol.~12, 2020.

\bibitem{BoZeng2020}
B.~Zeng, ``A practical scheme to compute the pessimistic bilevel optimization
  problem,'' {\em INFORMS Journal on Computing}, vol.~32, no.~4,
  pp.~1128--1142, 2020.

\bibitem{porta_1.4.1}
T.~Christof and A.~Löbel, ``Polyhedron representation transformation
  algorithm,'' 2015.
\newblock Software available from https://porta.zib.de/.

\bibitem{bejani2021systematic}
M.~M. Bejani and M.~Ghatee, ``A systematic review on overfitting control in
  shallow and deep neural networks,'' {\em Artificial Intelligence Review},
  vol.~54, no.~8, pp.~6391--6438, 2021.

\bibitem{montesinos2022overfitting}
O.~A. Montesinos~L{\'o}pez, A.~Montesinos~L{\'o}pez, and J.~Crossa,
  ``Overfitting, model tuning, and evaluation of prediction performance,'' in
  {\em Multivariate statistical machine learning methods for genomic
  prediction}, pp.~109--139, Springer, 2022.

\bibitem{bennouna2023}
A.~Bennouna and B.~V. Parys, ``Holistic robust data-driven decisions,'' {\em
  arXiv preprint arXiv:2207.09560}, 2023.

\end{thebibliography}

\end{document}